\documentclass[11pt, keyword, address]{article}
\usepackage{graphicx, comment}
\usepackage{times}
\usepackage{amsmath}
\usepackage{amssymb}
\usepackage{latexsym}
\usepackage{amsthm}
\newtheorem{theorem}{Theorem}[section]
\newtheorem{lemma}[theorem]{Lemma}
\newtheorem{proposition}[theorem]{Proposition}
 
\newtheorem{definition}[theorem]{Definition}
\newtheorem{corollary}[theorem]{Corollary}
\newtheorem{remark}[theorem]{Remark}
\newcommand{\As}{\mathcal{A}_{\mathrm{s}}}
\newcommand{\Aw}{\mathcal{A}_{\mathrm{w}}}
\newcommand{\A}{\mathcal{A}}
\newcommand{\Areg}{\mathcal{A}_{\mathrm{reg}}}
\newcommand{\ds}{\mathrm{d}_{\mathrm{s}}}
\newcommand{\dw}{\mathrm{d}_{\mathrm{w}}}
\newcommand{\dd}{\mathrm{d}}
\newcommand{\Dc}{\mathcal{E}}
\newcommand{\Fc}{\mathcal{F}}
\newcommand{\s}{\mathrm{s}}
\newcommand{\w}{\mathrm{w}}
\newcommand{\ww}{\omega_{\mathrm{w}}}
\newcommand{\ws}{\omega_{\mathrm{s}}}
\newcommand{\Xw}{X_{\mathrm{w}}}
\newcommand{\Hw}{H_{\mathrm{w}}}
\newcommand{\Ab}{\mathcal{A}_{\bullet}}
\newcommand{\db}{\mathrm{d}_{\bullet}}
\newcommand{\I}{\mathcal{I}}
\newcommand{\Bb}{B_{\bullet}}
\newcommand{\wb}{{\omega}_{\bullet}}
\newcommand{\ddt}{\frac{d}{dt}}

\newcommand{\ts}{\textstyle}
\def\setmarsing{
\oddsidemargin-0in
\evensidemargin-0in
\textwidth5.5in
\textheight8in
}
\setmarsing

\title{On global attractors of the 3D Navier-Stokes equations}
\author{A. Cheskidov\thanks{University of Michigan, Department of Mathematics,
Ann Arbor, MI 48109 (acheskid@umich.edu)} \ and
C. Foias\thanks{Texas A\&M University, Department of Mathematics,
College Station, TX 77843}}
\begin{document}
\maketitle
\abstract{
In view of the possibility that the 3D Navier-Stokes equations (NSE)  might
not always have regular solutions, we introduce an abstract
framework for studying the asymptotic behavior
of multi-valued dissipative evolutionary systems with respect to two
topologies---weak and strong. Each such system
possesses a global attractor in the weak topology, but not necessarily
in the strong. In case the latter exists and is weakly closed, it coincides
with the weak global attractor. We give a sufficient condition for
the existence of the strong global attractor, which is verified for the
3D NSE when all solutions on the weak global attractor are strongly continuous.
We also introduce and study a two-parameter family of models for the
Navier-Stokes equations, with similar properties and open problems.
These models always possess weak global attractors, but on some of them
every solution blows up (in a norm stronger than the standard energy one) in finite time.
}

\section{Introduction}
A remarkable feature of many dissipative partial differential equations (PDEs)
is the existence of a global attractor to which all the solutions converge as
time goes to infinity \cite{SY,T2}. The global attractor $\A$
is the minimal closed set in
a phase space $H$ (i.e., the functional space, usually a Banach space,
in which the solutions exist) that uniformly attracts the trajectories starting
from any a priori given
bounded set in $H$. When the topology on $H$ is referred as the strong
(weak) topology of $H$, we will call $\A$ the strong (respectively weak)
global attractor.

It is possible that a dissipative PDE does not have a strong global attractor.
For instance, the 2D Navier-Stokes equations (NSE) on a bounded domain
$\Omega \subset \mathbb{R}$, when supplemented with appropriate boundary
conditions, possess a strong global attractor in $H$ (a certain subspace of
$L^2(\Omega)^3$) \cite{L,T2}, but it is not yet known
whether this holds for the 3D NSE.

Nevertheless, even for the 3D NSE one can prove that there exists
a weak global attractor \cite{FT85}.
When the strong global attractor is strongly compact in $H$ (e.g., in the 2D NSE),
then it is also the weak global attractor. But, in any case, the weak global
attractor is an
appropriate generalization of the strong global attractor since it captures
the long-time behavior of the solutions. In particular, the support of any
time-average measure of the 3D NSE is included in the weak global attractor
(see \cite{FMRT}).
One should note that Sell \cite{S} introduced a related notion of a trajectory
attractor $\mathfrak{A}$ in the space of all trajectories, which was further
studied in \cite{CV, FS, S, SY}. The weak global attractor coincides with the
set of values of all trajectories in $\mathfrak{A}$ at any fixed time $t$.

The aim of this study is to present a general abstract framework which
is applicable to the 3D NSE even in the case where they do not
possess a strong global attractor. This framework may be also useful in
the study of other PDEs for which the existence of the strong global
attractor is in limbo. This aim forces us to consider multi-valued
evolutionary systems.
A number of papers have
been published concerning attractors of multi-valued semiflows. See
\cite{CMR} for a comparison of two canonical abstract frameworks by
Melnik and Valero \cite{MV} and Ball \cite{B1}. The main difference between
our evolutionary system $\Dc$ and Ball's generalized semiflow is that
we do not include
the hypotheses of concatenation and upper semicontinuity  with respect to
the initial data. This allows us to consider an evolutionary system whose
trajectories are all Leray-Hopf weak solutions of the 3D NSE.

Our definition of the evolutionary system $\Dc$ already exploits the effect of
dissipativity, namely the existence of an absorbing ball.
In fact, the space $X$ in which the trajectories of $\Dc$ live is, in applications,
precisely such an absorbing ball. Since in most applications the phase
space is a separable reflexive Banach space, both the strong and the
weak topologies on $X$ are metrizable. This is the motivation for us
to define strong and weak topologies on $X$ to be the ones induced by
appropriate metrics.

 We show
that every evolutionary system always possesses a weak global
attractor;  moreover, if the strong global attractor exists and is weakly closed,
then it has to coincide with the weak global attractor.
Note that some classical definitions (see, e.g., \cite{T2}) require a global
attractor to be an invariant set. We will see that under
a condition, which is, for example, satisfied by the Leray-Hopf weak solutions
of the 3D NSE, the weak global attractor is also the maximal bounded
invariant set. We recall that those solutions are always weakly continuous
in $L^2(\Omega)^3$.

It is known that if a weak global attractor for the 3D NSE is bounded in $V$,
then it is in fact strong \cite{FT85}. Moreover,
Ball \cite{B1} showed that if a generalized semiflow for a dissipative evolutionary
system is  asymptotically compact, then a strong global attractor exists.
This generalizes corresponding results for semiflows (see {\cite{H,HLS,L1}})
and implies that the strong global attractor for the 3D NSE exists
under the condition that all weak solutions are strongly continuous from $(0, \infty)$
to $L^2(\Omega)^3$ (see \cite{B1}). In this paper we show that even without the assumptions of concatenation and upper
semicontinuity  with respect to the initial data, the asymptotic compactness implies
that the weak global attractor is the minimal compact attracting set in the strong
metric, i.e., the weak global attractor is in fact the strong compact global attractor
(Theorem~\ref{t:asymptoticcompact}).
Applied to the 3D NSE, this result implies the existence of a strong compact global
attractor in the case when the solutions on the weak global attractor are
continuous in $L^2(\Omega)^3$ (Theorem~\ref{thm:main}).

The convergence of Leray-Hopf weak solutions was also studied by Rosa \cite{R}, 
namely, he introduced an asymptotic regularity condition
to insure the strong convergence of a weak solution towards its weak
$\omega$-limit. This condition requires the limit solutions to be strongly 
continuous in $L^2(\Omega)^3$, and implies that the weak global attractor is a
strongly compact strongly attracting set if the weak solutions on the weak global
attractor are strongly continuous in $L^2(\Omega)^3$.
Moreover, since a trajectory of the 3D NSE that is not strongly continuous
in  $L^2(\Omega)^3$ also, obviously, is not relatively strongly compact in $L^2(\Omega)^3$,
the strong continuity of weak solutions on the weak global
attractor $\Aw$ is also a necessary condition for $\Aw$ to be strongly compact
(see \cite{R}).

Recall that we define a global attractor as the minimal closed (uniformly)
attracting set in the corresponding metric, and hence, allow the possibility
for the solutions on a strong global attractor to be discontinuous.
We address this issue by
studying the weak convergence of weak solutions towards a weak solution strongly continuous from the right in $L^2(\Omega)^3$.  
We show that the weak convergence is strong under the
condition that the "energy jumps" of solutions converge to the "energy jumps"
of the limit solution (Theorem~\ref{t:deltat}). However, at this stage, no
necessary condition for the strong convergence is known.

Finally, we provide an example of a dissipative evolutionary system for
which all solutions on the weak global attractor blow-up in finite time.
We introduce a two-parameter family of simple infinite-dimensional
systems of differential equations.
These systems, called {\it tridiagonal models for the NSE (TNS models)},
display basic features of the NSE.
In particular, they are examples of dissipative systems that possess a weak
global attractor, but the existence of a strong global attractor is not
known.

TNS models have similar form and some similar properties to shell
models, specifically dyadic models studied in \cite{FP,KP}.
Moreover, a similar analysis of the dyadic models also results in a finite
time  blow-up (see \cite{C}).
In the TNS models though, the coefficients in the equations are chosen
to yield NSE-like scaling properties.
More precisely, we first mimic the Stokes operator in 3D via
the choice a positive definite operator on $l^2$ whose eigenvalues grow with the same speed as the eigenvalues of the Stokes operator. 
Second, we mimic the nonlinear term of the NSE via the choice of a bilinear
form, which scales like the Sobolev estimate for the NSE.
Then we obtain the following system of differential equations:
\begin{equation} \label{model-intro}
\ddt u + \nu Au + B(u,u) = g,
\end{equation}
where $u=(u_1,u_2,\dots)$,
\[
(Au)_n = n^{\alpha} u_n,
\]
and
\[
(B(u,v))_n=- n^\beta u_{n-1} v_{n-1} + (n+1)^\beta u_n v_{n+1},
\qquad n =1, 2, \dots,
\]
with $u_0=0$.
Here $\alpha$ and $\beta$ are two positive parameters. Note that
the orthogonality property in $l^2$ holds for $B$, which implies the existence
of an absorbing ball and a weak global attractor. Moreover, when
$\alpha = 2/3$, which corresponds to the speed with which the eigenvalues
of the Stokes operator grow in three-dimensional space, and $\beta=11/6$,
we have the following sharp estimate:
\[
| (B(u,u),Au) | \lesssim |Au |^{3/2} |A^{1/2}u|^{3/2},
\]
where
$
|v|^2=\sum v_n^2. 
$
This estimate is exactly the same as the the estimate based on the Sobolev
inequalities for the nonlinear term of the 3D NSE, with $|\cdot|$ being the $L^2$-norm.
It is an open question whether solutions of (\ref{model-intro})
can lose regularity in the case $(\alpha,\beta)=(2/3,11/6)$. However,
we show that in the nonviscous case $\nu=0$, for every $\alpha>0$, $\beta>0$,
$\gamma>0$,
and $g_n\geq 0$ for all $n\in\mathbb{N}$, the norm
$|A^{(\beta+\gamma-1)/(3\alpha)} u|$ of every solution
with $u_n(0)\geq 0$ and $u(0)\ne0$
blows up (Theorem~\ref{t:Eblow}).

When the viscosity is not zero, the model always possesses a weak global
attractor $\Aw$. Moreover, we prove that if the force $g$ is
large enough, then all solutions on $\Aw$  blow up in finite time in
an appropriate norm when $2\beta>3\alpha+3$ (Remark~\ref{rem:bl}).
The question whether $\Aw$ is a strong global attractor 
remains open in the case where $\beta > \alpha+1$.

\section{Evolutionary system and global attractors}

Let $(X,\ds(\cdot,\cdot))$ be a metric space endowed with
a metric $\ds$, which will be referred to as a strong metric.
Let $\dw(\cdot, \cdot)$ be another metric on $X$ satisfying
the following conditions:
\begin{enumerate}
\item $X$ is $\dw$-compact.
\item If $\ds(u_n, v_n) \to 0$ as $n \to \infty$ for some
$u_n, v_n \in X$, then $\dw(u_n, v_n) \to 0$ as $n \to \infty$,
that is, the identity map $(X, \ds) \mapsto (X,\dw)$ is uniformly continuous.
\end{enumerate}
Due to the latter property, $\dw$ will be referred to as the weak metric on $X$.
Note that any strongly compact ($\ds$-compact) set is weakly compact
($\dw$-compact), and any weakly closed set is strongly closed.
Also it will be convenient to denote by $\overline{A}^{\bullet}$ the
closure of the set $A\subset X$ in the topology generated by $\db$; here and
throughout $\bullet$ stands for either $\mathrm{s}$ or $\mathrm{w}$.

To define an evolutionary system, first let
\[
\mathcal{T} := \{ I: \ I=[T,\infty) \mbox{ for some } T\in \mathbb{R}, \mbox{ or } 
I=(-\infty, \infty) \},
\]
and for each $I \subset \mathcal{T}$ let $\mathcal{F}(I)$ denote
the set of all $X$-valued functions on $I$.
Now we define an evolutionary system $\Dc$ as follows.
\begin{definition} \label{d:E}
A map $\Dc$ that associates to each $I\in \mathcal{T}$ a subset
$\Dc(I) \subset \mathcal{F}$ will be called an evolutionary system if
the following conditions are satisfied:
\begin{enumerate}
\item $\Dc([0,\infty))\ne \emptyset$.
\item
$\Dc(I+s)=\{u(\cdot): \ u(\cdot -s) \in \Dc(I) \}$ for
all $s \in \mathbb{R}$.
\item $\{u(\cdot)|_{I_2} : u(\cdot) \in \Dc(I_1)\}
\subset \Dc(I_2)$ for all
pairs of $I_1,I_2 \in \mathcal{T}$, such that $I_2 \subset I_1$.
\item
$\Dc((-\infty , \infty)) = \{u(\cdot) : \ u(\cdot)|_{[T,\infty)}
\in \Dc([T, \infty)) \ \forall T \in \mathbb{R} \}.$
\end{enumerate}
\end{definition}
We will refer to $\Dc(I)$ as the set of all trajectories (solutions)
on the time interval $I$. Let $P(X)$ be the set of all subsets of $X$.
For every $t \geq 0$ define a map
\begin{eqnarray*}
&R(t):P(X) \to P(X),&\\
&R(t)K := \{u(t): u(0)\in K, u(\cdot) \in \Dc([0,\infty))\}, \qquad
K \subset X.&
\end{eqnarray*}
Note that the assumptions on $\Dc$ imply that $R(t)$ enjoys
the following property:
\begin{equation} \label{eq:propR(T)}
R(t+s)K \subset R(t)R(s)K, \qquad K \subset X,\quad t,s \geq 0.
\end{equation}

Let us first point out that the trajectories are not required to be continuous
(even with respect to $\dw$) nor are they uniquely determined by their
starting points, i.e.,  it is possible to have two different
trajectories $u,v \in \Dc([0,\infty))$, such that $u(0)=v(0)$. Second, there is no assumption of concatenation. If $u \in \Dc([0,\infty))$ and
$v \in \Dc([T,\infty))$ for some $T>0$, so that
$u(T)=v(T)$, then the following function
\[
w(t)= \left\{
\begin{aligned}
u(t), \qquad &\mbox{if } t\in[0,T],\\
v(t), \qquad &\mbox{if } t\in(T,\infty).
\end{aligned}
\right.
\]
need not be in $\Dc([0,\infty))$.
We avoid the assumptions of the continuity,
uniqueness, and concatenation in order to be able to consider
an evolutionary system consisting of Leray-Hopf weak solutions to the
3D Navier-Stokes equations.

Often an evolution of a dynamical system can be described by
a semigroup of continuous mappings acting on some metric space $H$:
\[
S(t): H \to H, \quad t \geq 0.
\]
The semigroup properties are the following:
\begin{equation} \label{eq:semigroup}
S(t+s) = S(t)S(s), \quad t,s \geq 0, \qquad S(0)= \mbox{Identity operator.}
\end{equation}
Then, if $u(t)\in H$ represents a state of the dynamical system at 
time $t$, we have
\[
u(t+s) = S(t)u(s), \qquad t,s \geq 0.
\]
A ball $B\subset H$ is called an absorbing ball, if for any bounded
set $K \subset H$ there exists $t_0$, such that
\[
S(t)K \subset B, \qquad \forall t \geq t_0.
\]
Hence, if we are interested in the long-time behavior of the dynamical system,
it is enough to consider a restriction of the system to an absorbing ball.
So, we let $X$ be a closed absorbing ball, and
define the map $\Dc$
in the following way:
\[
\Dc(I) := \{u(\cdot): u(t+s)=S(t)u(s) \mbox{ and }
u(s) \in X \ \forall s\in I, \ t \geq 0 \}.
\]
Note that conditions 1--4 for the evolutionary system $\Dc$ automatically follow
from the semigroup properties (\ref{eq:semigroup}) of $S(t)$.
In addition, let $T$ be such that
\[
S(t)X \subset X \qquad \forall t \geq T.
\]
Then we have
\[
R(t)K = S(t)K, \qquad \forall K \subset S(T)X, \ t \geq 0.
\]
The 3D Navier-Stokes equations will serve as an instructive illustration and 
application  of our consideration. 
As yet in our knlowledge of the 3D NSE, the time evolution of the 3D NSE cannot
be described by a semigroup of maps. Therefore, for the 3D NSE we will have a
more involved definition of $\Dc$ (see Section~\ref{3DNSE}).

Having defined the evolutionary system $\Dc$, we proceed to define attracting
sets and global attractors.
For $A \subset X$ and $r>0$, denote
$
B_{\bullet}(A,r) = \{u: \ \db(u, A) < r\},
$
where
\[
\db(u, A):=\inf_{x\in A}\db(u,x).
\]

\begin{definition}
A set $A \subset X$ is a $\mathrm{d}_{\bullet}$-attracting set
($\bullet = \mathrm{s,w}$) if it uniformly
attracts $X$ in $\mathrm{d}_{\bullet}$-metric, i.e.,
for any $\epsilon>0$, there exists
$t_0$, such that
\[
R(t)X \subset B_{\bullet}(A, \epsilon), \qquad \forall t \geq t_0.
\]
\end{definition}

\begin{definition}
$\mathcal{A}_{\bullet}\subset X$ is a
$\db$-global attractor ($\bullet = \mathrm{s,w}$) if
$\mathcal{A}_{\bullet}$ is a minimal $\db$-closed
$\db$-attracting  set, i.e., $\mathcal{A}_{\bullet}$
is $\db$-closed $\db$-attracting
and every subset $A \subset \Ab$ that is also $\db$-closed and $\db$-attracting
satisfies $A=\Ab$.
\end{definition}
Note that the empty set is never an attracting set.
Note also that since $X$ is not strongly compact, the intersection of two 
$\ds$-closed $\ds$-attracting sets might not be $\ds$-attracting.
Nevertheless, the uniqueness of a global attractor is a direct consequence
of the following lemma.

\begin{lemma}
If $\mathcal{A}_{\bullet}$ exists and $A$ is a $\db$-closed $\db$-attracting set,
then $\Ab \subset A$ ($\bullet = \mathrm{s,w}$).
\end{lemma}
\begin{proof}
Let $A$ be an arbitrary $\db$-closed $\db$-attracting set.
Take any point $a \in \Ab$. Let $\epsilon>0$. If there exists $t_\epsilon>0$, such that 
\[
R(t)X \cap \Bb(a,\epsilon) = \emptyset, \qquad \forall t\geq t_\epsilon,
\]
then $\Ab \setminus \Bb(a,\epsilon/2)$ is a $\db$-closed $\db$-attracting
set contradicting the minimality of $\Ab$. 
So, there exists a sequence $t_n \to \infty$ as $n\to \infty$,
such that
\[
R(t_n)X \cap \Bb(a,\epsilon) \ne \emptyset, \qquad \forall n.
\]
On the other hand, since $A$ is $\db$-attracting, we infer that
\[
R(t_n)X \subset \Bb(A,\epsilon),
\]
for $n$ large enough. It follows that
\[
A \cap \Bb(a,2\epsilon) \ne \emptyset.
\]
Since $A$ is $\db$-closed, we have that $a \in A$.  Thus,
$\Ab \subset A$.
\end{proof}
As a direct consequence of this lemma we have the following.

\begin{corollary}
If $\mathcal{A}_{\bullet}$ exists, then it is unique
($\bullet = \mathrm{s,w}$). 
\end{corollary}
Assume now that the weak global attractor $\Aw$ exists. If $\Aw$ is
a strongly attracting set, does it follow that the strong global attractor
exists? In general, this may not be true. However, if the strong global attractor
exists, then clearly it is a $\dw$-attracting set. Moreover, we have the following.

\begin{theorem} \label{thm:coinside}
If $\As$ exists, then $\Aw$ exists and
\[
\Aw=\overline{\As}^{\w}.
\]
\end{theorem}
\begin{proof}
If there exists
a $\dw$-closed $\dw$-attracting set $A\subset \overline{\As}^{\mathrm{w}}$ and
$A\ne \overline{\As}^{\mathrm{w}}$, then there exists
$u_0 \in \As$, such that
\[
d=\dw(u_0, A)>0.
\]
By the definition of an attracting set, there exists a time $t_0>0$, such that
\begin{equation} \label{e:111}
R(t)X \subset B_{\mathrm{w}}(A, d/2) \qquad \forall t\geq t_0.
\end{equation}
Note that
\[
\dw(u_0, B_{\mathrm{w}}(A, d/2)) \geq d/2.
\]
Therefore, by virtue of Property 2 in the definition of $\dw$, there
exists $\delta>0$, such that
\[
\ds(u_0, B_{\mathrm{w}}(A, d/2)) > \delta,
\]
whence,
\[
B_{\mathrm{s}}(u_0, \delta) \cap B_{\mathrm{w}}(A, d/2) = \emptyset.
\]
Now from \eqref{e:111} it follows that
\[
B_{\mathrm{s}}(u_0, \delta)\cap R(t)X =\emptyset \qquad \forall t \geq t_0.
\]
Consequently, $\As \setminus B_{\mathrm{s}}(u_0, \delta/2)$ is a $\ds$-closed 
$\ds$-attracting set strictly included in $\As$, a contradiction.
Hence, $\overline{\As}^{\mathrm{w}}$ is the weak global attractor.
\end{proof}

The following are two simple examples of evolutionary systems
that possess a weak global attractor $\Aw$, but not a strong global
attractor $\As$.

\noindent
\textbf{Example 1.}
Let
\[
X=\left\{ u \in L^2(-\infty,\infty): \int_{-\infty}^\infty u(x)^2 \, dx \leq 1,
u(x)=0 \mbox{ for } x> 0\right\},
\]
and define on $X$ the distances
\[
\ds(u,v):=\left( \int_{-\infty}^\infty(u(x)-v(x))^2 \, dx\right)^{1/2}, \qquad
\dw(u,v)=\int_{-\infty}^\infty \frac{1}{2^{|x|}}
\frac{|u(x)-v(x)|}{1 + |u(x)-v(x)|} \, dx.
\]
Consider the following partial differential equation:
\[
\frac{\partial u}{\partial t} = \frac{\partial u}{\partial x}.
\]
The trajectories of the evolutionary system $\Dc$
will be solutions of this equation, i.e.,
\[
\begin{aligned}
&\Dc([s,\infty)) =\{u\in \Fc([s,\infty)) : u(t) =u_0(\cdot+t-s), t\in[s,\infty), u_0 \in X\}, \qquad
\forall s \in \mathbb{R},\\
& \Dc ((-\infty,\infty))=\{0\}.
\end{aligned}
\]
Then it is easy to show that $\Aw=\{0\}$ (see also Theorem~\ref{c:weakA}). However, no trajectory except the trivial one $u=0$ strongly
converges to $0$ as $t\to \infty$.

\noindent
\textbf{Example 2.}
Take
\[
X=\left\{ u\in l^2: \sum_{n=1}^\infty u_n^2  \leq 1\right\}.
\]
where $u=(u_n)$, $u_n \in \mathbb{R}$ for all $n$.  For $u, v \in X$, let
\[
\ds(u,v):=\left( \sum_{n=1}^\infty(u_n-v_n)^2 \right)^{1/2}, \qquad
\dw(u,v)= \sum_{n=1}^\infty \frac{1}{2^n}
\frac{|u_n-v_n|}{1 + |u_n-v_n|}.
\]
Consider the following differential equation:
\[
\ddt u_n = -\frac{1}{n} u_n, \qquad n \in \mathbb{N}
\]
The trajectories of the evolutionary system $\Dc$
will be solutions of this equation, i.e.,
\[
\Dc([s,\infty)) =\{u\in\Fc([s,\infty)): (u_n(t))=(u_n^0e^{(s-t)/n}), t\in[s,\infty),
u^0 \in X\},
\]
for $s \in \mathbb{R}$, and
\[
\Dc((-\infty,\infty))=\{0\}.
\]
Take any $u^0 \in X$ and consider the trajectory $u \in \Dc([0,\infty))$ starting
at $u^0$, i.e., $u(0)=u^0$. Then we have
\[
\ds(u(t),0)^2=\sum_{n=1}^\infty u_n(t)^2 = \sum_{n=1}^\infty(u^0_n)^2 e^{-\frac{2t}{n}}
\to 0 \qquad \mbox{as} \qquad t \to \infty.
\]
However, the convergence is not uniform in the $\ds$-metric,
although it is uniform in the $\dw$-metric. So
again $\Aw=\{0\}$, but $\As$ does not exist.

Note that the nonexistence of $\As$ in the two examples is due to
two different behaviors of the trajectories.
In the first example all the nontrivial trajectories converge
to $\Aw$ weakly, but not strongly. In the second example all the nontrivial
trajectories converge to $\Aw$ strongly, but not uniformly.

\begin{definition} \label{d:ucomp}
The map $R(t)$ is uniformly $\db$-compact
($\bullet= \mathrm{s, w}$) if
there exists $t_0\geq 0$, such that
\[
\bigcup_{t \geq t_0} R(t) X
\]
is relatively $\db$-compact.
\end{definition}

Note that since $X$ is $\dw$-compact, 
$R(t)$ is automatically uniformly $\dw$-compact.

\begin{definition} The $\wb$-limit ($\bullet= \mathrm{s, w}$)
of a set $K \subset X$ is
\[
\wb (K):=\bigcap_{T\geq 0} \overline{\bigcup_{t \geq T} R(t) K}^{\bullet},
\]
where the closure is taken in $\db$-metric.
\end{definition}

\begin{lemma} \label{l:wbinA}
Let $A$ be a $\db$-closed $\db$-attracting set. Then
\[
\wb(X) \subset A.
\]
\end{lemma}
\begin{proof}
Suppose that there exists $a \in \wb(X) \setminus A$.
Since $A$ is $\db$-closed, there exists $\epsilon>0$, such that
\[
A \cap \Bb(a,\epsilon) = \emptyset.
\]
By the definition of the $\wb$-limit,
there exist a sequence $t_n \to \infty$ as $n \to \infty$ and
a sequence $x_n \in R(t_n)X$, such that $\db(x_n,a) \to 0$
as $n \to \infty$. Hence, there exists $N>0$, such that 
\[
x_n \notin \Bb(A, \epsilon/2), \qquad \forall n\geq N.
\]
This means that $A$ is not $\db$-attracting, a contradiction.
\end{proof}

\begin{lemma} \label{l:oncompactsemigroup}
If the map $R(t)$ is uniformly $\db$-compact
($\bullet= \mathrm{s, w}$), then 
$\wb(X)$ is a nonempty $\db$-compact $\db$-attracting set.
\end{lemma}
\begin{proof}
By Definition~\ref{d:ucomp} and the fact that $\Dc([0,\infty)) \ne \emptyset$,
there exists $t_0$, such that
\[
W(T):=\overline{\bigcup_{t \geq T} R(t) X}^{\bullet}
\]
is a nonempty $\db$-compact set for all $T\geq t_0$.
In addition, $W(s) \subset W(t)$ for all $s\geq t \geq 0$. Thus,
\[
\wb(X)= \bigcap_{T\geq t_0} W(T)
\]
is a nonempty  $\db$-compact set.

We will now prove that $\wb (X)$ uniformly $\db$-attracts $X$. Assume it
does not. Then there exists $\epsilon>0$, such that
\[
V(t):=W(t)\cap(X\setminus \Bb(\wb(X), \epsilon)) \ne \emptyset, \qquad \forall t\geq 0.
\]
Since $V(t)$ is $\db$-compact for $t\geq t_0$ and $V(s) \subset V(t)$ for
$s\geq t \geq 0$, we have that there exists
\[
x \in \bigcap_{t\geq t_0} V(t).
\]
Hence, $x \in \wb(X)$. However, this implies that $x\notin V(t)$, $t\geq 0$,
a contradiction.
\end{proof}

\begin{theorem} \label{thm:exofA}
If the map $R(t)$ is uniformly $\db$-compact
($\bullet= \mathrm{s, w}$), then the  $\db$-global
attractor exists and satisfies the following additional properties:
\begin{enumerate}
\item[(a)]
$\Ab = \wb(X)$.
\item[(b)]
$\Ab$ is $\db$-compact.
\end{enumerate}
\end{theorem}
\begin{proof}
By Lemma~\ref{l:oncompactsemigroup},
$\wb(X)$ is a nonempty $\db$-compact $\db$-attracting set.
Moreover, $\wb(X)$ is the
minimal $\db$-closed $\db$-attracting set due to Lemma~\ref{l:wbinA}.
Therefore, $\wb(X)$ is the $\db$-global attractor.
\end{proof}

Note that since $X$ is $\dw$-compact (see the definition of $\dw$), $R(t)$ is uniformly weakly compact.
Hence we have the following.
\begin{corollary}
The evolutionary system $\Dc$ always possesses a weak global attractor $\Aw$.
\end{corollary}
Our next goal is to investigate whether $\Aw$ is an invariant set in the following
sense.

\begin{definition} The set $A \subset X$ is invariant, if
\[
\{u(t): u\in \Dc((-\infty,\infty)), u(0) \in A \}= A, \qquad \forall t\geq 0.
\]
\end{definition}
Assume that $u_0$ belongs to some invariant set. Then for all $t>0$
we have that $u_0 \in R(t)X$.
Hence, $u_0 \in \omega_{\mathrm{w}}(X)=\Aw$. Therefore, $\Aw$ contains every invariant set. Moreover, we will show that $\Aw$ is invariant under some compactness
property that is for instance satisfied by the family of all 
Leray-Hopf solutions of the 3D NSE (see Section~\ref{3DNSE}).

Let $C([a, b];X_\bullet)$ be the space of $\db$-continuous $X$-valued
functions on $[a, b]$ endowed with the metric
\[
\dd_{C([a, b]; X_\bullet)}(u,v) = \sup_{t\in[a,b]}\db(u(t),v(t)). 
\]
Let also $C([a, \infty);X_\bullet)$ be the space of $\db$-continuous
$X$-valued functions on $[a, \infty)$
endowed with the metric
\[
\dd_{C([a, \infty); X_\bullet)}(u,v) = \sum_{T\in \mathbb{N}} \frac{1}{2^T} \frac{\sup\{\db(u(t),v(t)):a\leq t\leq a+T\}}
{1+\sup\{\db(u(t),v(t)):a\leq t\leq a+T\}}.
\]

\begin{theorem} \label{c:weakA}
If $\Dc([0,\infty))$ is compact in $C([0, \infty);\Xw)$, then
\begin{enumerate}
\item[(a)]
$
\Aw = \mathcal{I}:=\{ u_0: \ u_0=u(0) \mbox{ for some }
u \in \Dc((-\infty, \infty))\}.
$
\item[(b)] $\Aw$ is the maximal invariant set.
\end{enumerate}
\end{theorem}
\begin{proof}
Since obviously $\I$ is the maximal invariant set, we have that $\I \subset \Aw$.
It remains to prove that $\Aw \subset \I$. Take any $a \in \Aw$.
Since $\Aw = \omega_{\mathrm{w}}(X)$,  there exist
$t_n \to \infty$, as $n\to \infty$ and $a_n\in R(t_n)X$, such that
$a_n \to a$ weakly as $n \to \infty$. 
Using Property 2 in Definition~\ref{d:E}, there exist
$u_n \in \Dc([-t_n, \infty))$, such that
$u_n(0) =a_n$. Also, Properties 2 and 3 in Definition~\ref{d:E}
of $\Dc$ imply that
$\Dc([-t_n,\infty))$ is compact in $C([-t_n, \infty);\Xw)$
and
\[
\{u|_{[-t_1,\infty)} : u \in \Dc([-t_{n},\infty))\}
\subset \Dc([-t_1,\infty)),
\]
for every $n$.
Now, passing to a subsequence and dropping a subindex, we can assume
that $u_n|_{[-t_1,\infty)} \to u^1\in \Dc([-t_1, \infty))$ in  $C([-t_1, \infty);\Xw)$
as $n\to \infty$. 
By a standard diagonalization process we obtain that there exist
$u \in  \Fc((-\infty,\infty))$ and a subsequence of $u_n$, still denoted by
$u_n$, such that 
$u_n|_{[-T,\infty)} \to u_{[-T,\infty)}$ in $C([-T, \infty);\Xw)$ for all $T>0$.
Thus, by the compactness we have that $u|_{[-T,\infty)} \in \Dc([-T, \infty))$
for all $T>0$, and hence $u \in \Dc((-\infty, \infty))$.
Finally, since $u(0)=a$. we have $a \in \mathcal{I}$.
Hence, $\Aw\subset \I$.
\end{proof}

Now that we know that the weak global attractor always exists, we can
weaken the condition on the existence of the strong global attractor. 

\begin{definition}
The map $R(t)$ is asymptotically $\db$-compact
($\bullet= \mathrm{s, w}$) if for any $t_n \to \infty$ and
any $x_n \in R(t_n) X$, the  sequence
$\{x_n\}$ is relatively $\db$-compact.
\end{definition}

\begin{theorem} \label{t:asymptoticcompact}
If the map $R(t)$ is asymptotically  $\ds$-compact, then
$\Aw$ is $\ds$-compact strong global attractor.
\end{theorem}
\begin{proof}
First note that $\ws(X)\subset \ww(X)=\Aw$. On the other hand,
let $a \in \Aw=\ww(X)$. By
the definition of $\ww$-limit, there exist $t_n \to \infty$ as $n \to \infty$ and
$x_n \in R(t_n)X$, such that 
\[
\dw(x_n, a) \to 0 \qquad \mbox{ as } \qquad n \to \infty.
\]
Thanks to the asymptotic compactness of $R(t)$, this convergence
is in fact strong. Therefore, $a\in \ws(X)$. Hence, $\ws(X)= \Aw$.

Now let us show that $\ws(X)$ is a $\ds$-attracting set. Assume that it is not. Then there
exist $\epsilon>0$, $x_n \in X$, and $t_n \to \infty$ as
$n\to \infty$, such that
\[
x_n \in R(t_n)X \setminus B_{\s}(\ws(X), \epsilon), \qquad
\forall n \in \mathbb{N}.
\]
Since  $R(t)$ is asymptotically $\ds$-compact, then
$\{x_n\}$ is relatively $\ds$-compact.
Passing to a subsequence and dropping a subindex, we may assume that
\[
x_n \to x \in X \qquad \mbox{strongly},
 \qquad \mbox{as} \ n \to \infty.
\]
Therefore, we have that $x  \in \omega_{\mathrm{s}}(X)$,
a contradiction.

Now note that $\ws(X)$ is the
minimal $\ds$-closed $\ds$-attracting set due to Lemma~\ref{l:wbinA}.
Therefore, $\ws(X)$ is the strong global attractor $\As$.
Finally, let us show that $\ws(X)$
is strongly compact. Take any sequence $a_n \in \ws(X)$. By the definition of
$\ws$-limit, there exist $t_n \to \infty$ and
$x_n \in R(t_n)X$, such that 
\[
\ds(x_n, a_n) \to 0 \qquad \mbox{ as } \qquad n \to \infty.
\]
Note that $\{x_n\}$ is relatively $\ds$-compact due to the asymptotic compactness
of $R(t)$. Hence, $\{a_n\}$ is relatively $\ds$-compact, which concludes the proof.

\end{proof}

Finally, in the following example we show that the asymptotic compactness
is a weaker condition than the uniform strong compactness.

\noindent
\textbf{Example.}
Take
\[
X=\left\{ u\in l^2: \sum_{n=1}^\infty u_n^2  \leq 1\right\}.
\]
where $u=(u_n)$, $u_n \in \mathbb{R}$ for all $n$.  For $u, v \in X$, let
\[
\ds(u,v):=\left( \sum_{n=1}^\infty(u_n-v_n)^2 \right)^{1/2}, \qquad
\dw(u,v)= \sum_{n=1}^\infty \frac{1}{2^n}
\frac{|u_n-v_n|}{1 + |u_n-v_n|}.
\]
Consider the following differential equations:
\[
\ddt u_1=0,
\]
and
\[
\ddt u_n = - u_n, \qquad n \in \mathbb{N}.
\]
The trajectories of the evolutionary system $\Dc$ will be solutions of this
equation, i.e.,
\[
\begin{split}
\Dc([s,\infty)) =\{&u\in\Fc([s,\infty)): u_1(t)=u^0_1, u_n(t)=u_n^0e^{(s-t)}
\mbox{ for } n\geq 2,\\ & t\in[s,\infty), u^0 \in X\},
\end{split}
\]
for $s \in \mathbb{R}$, and
\[
\Dc((-\infty,\infty))=\{u: u_1\in [-1,1], u_n=0, n\geq 2\}.
\]
Clearly, $R(t)$ is not uniformly strongly compact, but it is
asymptotically compact. Hence, the strong global attractor
exists and
\[
\As = \{u: u_1\in [-1,1], u_n=0, n\geq 2\}.
\]

Finally,
Theorems~\ref{c:weakA} and \ref{t:asymptoticcompact}
imply the following
\begin{remark} 
If $\Dc([0,\infty))$ is compact in $C([0, \infty);\Xw)$ and $R(t)$ is asymptotically  
$\ds$-compact, then  the strong global attractor $\As$ exists, and is the strongly compact
maximal invariant set. Consequently, $\As$ is a compact global attractor in the
conventional sense.
\end{remark}

\section{3D Navier-Stokes equations} \label{3DNSE}

Here we apply the results from the previous section 
to the space periodic 3D Navier-Stokes equations (NSE)
\begin{equation} \label{NSE1}
\left\{
\begin{aligned}
&\ddt u - \nu \Delta u + (u \cdot \nabla)u + \nabla p = f,\\
&\nabla \cdot u =0,\\
&u, p, f \mbox{ are periodic with period } L \mbox{ in each space variable,}\\
&u, f \mbox{ are in } L^2_{\mathrm{loc}}(\mathbb{R}^3)^3,\\
&u|_{t=0}=u_0,
\end{aligned}
\right.
\end{equation}
where $u$, the velocity, and $p$, the pressure, are unknowns, $f$ is
a given driving force, and $\nu>0$ is the kinematic  viscosity coefficient
of the fluid. By a Galilean change of
variables, we can assume that the space average of $u$ is
zero, i.e.,
\[
\int_\Omega u(x,t) \, dx =0, \qquad \forall t,
\]
where $\Omega=[0,L]^3$ is a periodic box.

In this section we will apply our general results from the previous section
to the study of the assymptotical behaviour of weak solutions to \eqref{NSE1}.

\subsection{Functional setting} \label{funcset}

First, let us introduce some notations and functional setting for \eqref{NSE1}.
We denote by
$(\cdot,\cdot)$ and $|\cdot|$ the $L^2(\Omega)^3$-inner product and the
corresponding $L^2(\Omega)^3$-norm.
Let $\mathcal{V}$ be the space of all $\mathbb{R}^3$ trigonometric polynomials of
period $L$ in each variable satisfying 
$\nabla \cdot u =0$ and $\int_\Omega u(x) \, dx =0$.
Let $H$ and $V$ to be the
closures of $\mathcal{V}$ in $L^2(\Omega)^3$ and $H^1(\Omega)^3$, respectively.
Also, define the distances $\db$ by
\[
\ds(u,v):=|u-v|, \qquad
\dw(u,v)= \sum_{\kappa \in \mathbb{Z}^3} \frac{1}{2^{|\kappa|}}
\frac{|u_{\kappa}-v_{\kappa}|}{1 + |u_{\kappa}-v_{\kappa}|},
\qquad u,v \in H,
\]
where $u_{\kappa}$ and $v_{\kappa}$ are Fourier coefficients of $u$
and $v$ respectively.

We denote by $P_{\sigma} : L^2(\Omega)^3 \to H$ the $L^2$-orthogonal
projection, referred to as the Leray projector, and by
$A=-P_{\sigma}\Delta = -\Delta$ the Stokes operator with the domain
$D(A)=(H^2(\Omega))^3 \cap V$. The Stokes operator is a self-adjoint
positive operator with a compact inverse.
Denote
\[
\|u\| := |A^{1/2} u| = \left(\int_\Omega \sum_{i,j=1}^3\left|\frac{\partial u_i}{\partial x_j}\right|^2 \, dx \right)^{1/2}.
\]
Note that $\|u\|$ is equivalent to the $H^1$-norm of $u$ for $u\in D(A^{1/2})$.

For a rigorous mathematical study of  the equation \eqref{NSE1} we need
a few more concepts from functional analysis. Namely, let $V'$ be the set
of all distributions of the form $v=\Delta u$, with $u\in V$. The $V'$-norm of
this $v$ is by definition $\|u\|$. Endowed with this norm, $V'$ becomes
the dual space of $V$, and if $v\in H$, then the value of $v$ at a point $w\in H$
equals to the usual scalar product $(v,w)$ in $H$. Now for $u$ and $v$ in $V$,
let $B(u,v):=P_{\sigma}(u \cdot \nabla v)$, which is an element of $V'$. If
$v \in \mathcal{D}(A)$, then $B(u,v) \in H$. Moreover, 
\[
\langle B(u,v),w \rangle =- \langle B(u,w),v \rangle, \qquad u,v,w \in V,
\]
in particular, $\langle B(u,v),v \rangle=0$ for all $u,v \in V$.

Equations (\ref{NSE1}) now can be condensed in the functional
differential equation
\begin{equation} \label{NSE}
\ddt u + \nu Au +B(u,u) = g \qquad \mbox{in} \qquad  V',
\end{equation}
where $u$ is a $V$-valued function of time and $g = P_{\sigma} f$.
Throughout, we will assume that  $g$ is time independent and $g\in H$.

\begin{definition}
A weak solution  of \eqref{NSE} on $[T,\infty)$ is an $H$-valued
function $u(t)$ defined for $t \in [T, \infty)$, such that
\[
u(\cdot) \in C([T, \infty); \Hw) \cap 
L_{\mathrm{loc}}^2([T, \infty); V),
\]
and
\begin{equation} \label{333}
w(\cdot):=g-\nu A u(\cdot) - B(u(\cdot),u(\cdot))\in L_{\mathrm{loc}}^1([T, \infty); V'),
\end{equation}
\begin{equation} \label{444}
\begin{split}
(u(t)-u(T),v)&= \langle u(t)-u(T),v \rangle\\
&=\int_T^t\langle w(s),v\rangle \, ds \qquad
\forall v \in V, t\geq T.
\end{split}
\end{equation}
\end{definition}
The relations \eqref{333} and \eqref{444} imply that
$\ddt u$ exists in $V'$ a.e. in $[T,\infty)$. Therefore, often in
the literature \eqref{444} and \eqref{333} are written as \eqref{NSE} and
\[
\ddt u \in L_{\mathrm{loc}}^1([T, \infty); V'),
\] 
respectively.

The classical fundamental result concerning \eqref{NSE} is the following.

\begin{theorem}[Leray, Hopf] \label{thm:Leray}
For every $u_0 \in H$,
there exists a weak solution $u(t)$ of (\ref{NSE}) on $[T,\infty)$ with $u(T)=u_0$
satisfying the following energy inequality: 
\begin{equation} \label{EI}
|u(t)|^2 + 2\nu \int_{t_0}^t \|u(s)\|^2 \, ds \leq
|u(t_0)|^2 + 2\int_{t_0}^t (g(s), u(s)) \, ds
\end{equation}
for all $t \geq t_0$, $t_0$ a.e. in $[T,\infty)$.
\end{theorem}
See \cite{K} for a hypothesis under which the energy equality holds. However,
in general, the existence of weak solutions satisfying the energy equality is
not known.  Therefore, we introduce the following definition. 
\begin{definition} \label{d:ex}
A Leray-Hopf solution of the \eqref{NSE} on the interval $[T, \infty)$
is a weak solution on $[T,\infty)$ satisfying the
energy inequality (\ref{EI}) for all $T \leq t_0 \leq t$,
$t_0$ a.e. in $[T,\infty)$. The set $Ex$ of measure $0$ of points $t_0$ for
which the energy
inequality does not hold will be called the exceptional set (of the solution).
In addition, the solution $u(t)$ will be called regular on an interval $(\alpha, \beta) \subset [T,\infty)$ if $u(t)\in V$ and $\|u(t)\|$ is continuous on $(\alpha,\beta)$.
\end{definition}

Note that the uniqueness of Leray-Hopf solutions of the Initial Value Problem is
not known.

\begin{theorem}[Leray] \label{existence}
For every $u_0 \in V$, there exists a strong solution
$u(t)$ of \eqref{NSE} on some interval $[0,T)$, $T>0$, with $u(0)=u_0$.
\end{theorem}

\begin{theorem}[Leray] \label{structure}
Let $u(t)$ be a Leray-Hopf solution of \eqref{NSE} on $[T,\infty)$. Then there are
at most countably many distinct open intervals $I_j$, such that
\[
[T,\infty) = \bigcup_j \overline{I}_j,
\]
$u(t)$ is regular on each $I_j$, and the measure of $[T,\infty)\setminus \cup_j I_j$
is zero.
\end{theorem}

\begin{theorem} \label{uniqueness}
Let $u(t)$ be a regular solution on $[0,T)$. Then every Leray-Hopf solution
$v(t)$ on $[0,\infty)$ with $v(0)=u(0)$ coincides with $u(t)$ in $[0,T)$.
\end{theorem}

Finally, we recall several well-known supplementary facts.
\begin{remark} \label{r:ex}
The complement of the exceptional set $Ex$ coincides with the set of points of
strong continuity from the right.
\end{remark}

\begin{theorem} \label{t:doubleconv}
Let $u(t)$, $u_n(t)$ be Leray-Hopf solutions on the interval $[0,\infty)$,
such that 
\[
u_n \to u \qquad \mbox{in} \qquad C([0,T];\Hw),
\]
as $n\to \infty$, for some $T>0$. Let $(\alpha,\beta) \subset (0, T)$ be an interval of regularity
of $u(t)$. Then for every $0 <\delta <(\beta-\alpha)/2$,
\[
\|u_n(t)-u(t)\|\to 0 \qquad \mbox{uniformly on} \qquad [\alpha +\delta, \beta -\delta],
\]
 as $n \to \infty$.
\end{theorem}

\begin{remark}
Let $u(t)$ be a Leray-Hopf solution. As a $V$-valued function, $u(t)$ is analytic in time on every interval of regularity.
\end{remark}
\subsection{The weak global attractor for the 3D NSE} \label{s:weakattractor}
A ball $B\subset H$ is called an absorbing ball for the equation \eqref{NSE}
if for any bounded set $K \subset H$, there exists $t_0$, such that
\[
u(t) \in B, \qquad \forall t \geq t_0,
\]
for all Leray-Hopf solutions $u(t)$ of \eqref{NSE} on $[0,\infty)$ with $u(0)\in K$.
It is well known that there exists an absorbing ball in $H$ for the 3D NSE.
In fact, one has the following.
\begin{proposition} \label{p:aset}
The 3D Navier-Stokes equations possess an absorbing ball
\[
B = B_{\mathrm{s}}(0, R),
\]
where $R$ is any number larger than
$|g|\nu^{-1} L/(2\pi)$ (see, e.g., \cite{CF}).
\end{proposition}

Fix $R>|g|\nu^{-1} L/(2\pi)$ and let $X$ be the closed absorbing ball
\[
X= \{u\in H: |u| \leq R\},
\]
which is, clearly, weakly compact. Then for any bounded set $K \subset H$,
there exists a time $t_0$, such that
\[
u(t) \in X, \qquad \forall t\geq t_0,
\]
for every Leray-Hopf solution $u(t)$ with the initial data $u(0) \in K$.
Classical NSE estimates (see \cite{CF}) imply that for any sequence of Leray--Hopf solutions $u_n(t)$ (not only for the ones guaranteed by Theorem~\ref{thm:Leray}),
the following result holds.

\begin{lemma} \label{l:convergenceofLH}
Let $u_n(t)$ be a sequence of Leray-Hopf solutions,
such that $u_n(t) \in X$ for all $t\geq t_0$. Then 
\[
\begin{aligned}
u_n \ \ &\mbox{is bounded in} \ \ L^2([t_0,T];V),\\
\ddt u_n \ \  &\mbox{is bounded in} \ \ L^{4/3}([t_0,T];V'),
\end{aligned}
\]
for all $T>t_0$.
Moreover, there exists a subsequence $u_{n_j}$ of $u_n$, which converges
in $C([t_0, T]; \Hw)$ to some Leray-Hopf solution $u(t)$, i.e.,
\[
(u_{n_j},v) \to (u,v) \qquad
\mbox{uniformly on} \qquad  [t_0,T],
\]
as $n_j\to \infty$, for all $v \in H$.
\end{lemma}

Consider an evolutionary system for which
a family of all trajectories consists
of all Leray-Hopf solutions of the 3D Navier-Stokes equations
in $X$. More precisely, define
\[
\begin{split}
\Dc([T,\infty)) := \{&u(\cdot): u(\cdot)
\mbox{ is a Leray-Hopf solution on } [T,\infty)\\
& \mbox{and } u(t) \in X \ \forall t \in [T,\infty)\},
\qquad T \in \mathbb{R},
\end{split}
\]
\[
\begin{split}
\Dc((\infty,\infty)) :=\{&u(\cdot): u(\cdot)
\mbox{ is a Leray-Hopf solution on } (-\infty,\infty)\\
& \mbox{and } u(t) \in X \ \forall t \in (-\infty,\infty)\}.
\end{split}
\]
Since $X$ is weakly compact, the existence of the weak
global attractor is a direct consequence of Theorem~\ref{c:weakA}.
Moreover, we have the following.
\begin{lemma} \label{l:compact}
$\Dc([0,\infty))$ is compact in $C([0, \infty);\Hw)$.
\end{lemma}
\begin{proof}
Take any sequence
$u_n \in \Dc([0,\infty))$, $n\in \mathbb{N}$.
Thanks to Lemma~\ref{l:convergenceofLH}, there exists
a subsequence, still denoted by  $u_n$, that converges
to some $u^{1} \in \Dc([0,\infty))$ in $C([0, 1];\Hw)$ as $n \to \infty$.
Now, passing to a subsequence and dropping a subindex, we obtain that
there exists $u^{2}\in \Dc([0,\infty))$, such that
$u_n \to u^2$ in $C([0, 2];\Hw)$ as $n \to \infty$.
Note that $u^1(t)=u^2(t)$ on $[0, 1]$.
Continuing
this diagonalization process, we obtain a subsequence $u_{n_j}$
of $u_n$ that converges
to some $u \in \Dc([0,\infty))$ in $C([0, \infty);\Hw)$ as $n_j \to \infty$, which
concludes the proof.
\end{proof}

Now Theorem~\ref{c:weakA} yields the following.
\begin{theorem} \label{t:AweakNSE}
The weak global attractor $\Aw$ for the 3D
Navier-Stokes equations exists and satisfies
\begin{enumerate}
\item[(a)]
$
\Aw = \{u(0): u \in \Dc((-\infty, \infty))\}.
$
\item[(b)] $\Aw$ is the maximal invariant set.
\end{enumerate}
\end{theorem}

\begin{lemma}
If $u(t)$, a Leray-Hopf solution of the 3D NSE, satisfies
\[
\limsup_{t \to \infty} \|u(t)\| < \infty,
\]
then $u(t)$ converges strongly in $H$ to the weak global
attractor $\Aw$.
\end{lemma}
\begin{proof}
Suppose that $u(t)$ does not converge strongly in $H$ to $\Aw$. 
Then there exist $M>0$ and a sequence $t_n \to \infty$ as
$n \to \infty$, such that
\begin{equation} \label{eq:farfromAw}
\ds(u(t_n), \Aw) > M, \qquad n \in \mathbb{N}.
\end{equation} 
Note that there exists a time $T>0$, such that
$\{u(t): t\geq T\}$ is relatively compact in $H$. Therefore, passing
to a subsequence, we may assume that $u(t_n)$ converges strongly
(and hence weakly) in $H$ to some $a \in H$. Therefore, $a \in \Aw$, which
contradicts
(\ref{eq:farfromAw}).
\end{proof}

\subsection{Strong convergence of Leray-Hopf solutions}

The aim of this subsection is to give sufficient conditions for a
sequence of Leray-Hopf solutions on $[T,\infty)$ to converge
in $C([T,\infty);H)$, provided it converges in $C([T,\infty);\Hw)$.
We start with preliminary properties, some of which may have
an intrinsic interest.
\begin{theorem} \label{t:lim}
Let $u(t)$ be a Leray-Hopf solution of \eqref{NSE} on $[T,\infty)$.
Let $Ex$ be the exceptional set for this
solution (see Def.~\ref{d:ex}). Then for any time $t_0 >T$,
there exist $A_-$ and $A_+$, such that
\begin{enumerate}
\item[(a)]
For every sequence $\{t_n\} \subset [T,\infty) \setminus Ex$,
such that $t_n \to t_0$, $t_n < t_0$, it follows that  $|u(t_n)| \to A_-$
as $n \to \infty$.
\item[(b)]
For every sequence $\{t_n\} \subset [T,\infty) \setminus Ex$,
such that $t_n \to t_0$, $t_n > t_0$, it follows that $|u(t_n)| \to A_+$
as $n \to \infty$.
\end{enumerate}
For these $A_-$ and $A_+$ we will use the following notations:
\[
\widetilde{\lim_{t \to t_0 -}} |u(t)| := A_- \qquad \mbox{and}
\qquad \widetilde{\lim_{t \to t_0 +}} |u(t)| := A_+.
\]
\end{theorem}
\begin{proof}
For $\{t_n\}$ as in (a), the energy inequality (\ref{EI}) on
$[t_n, t_{n+k}]$ is 
\[
|u(t_{n+k})|^2 + 2\nu \int_{t_n}^{t_{n+k}}\|u(s)\|^2 \, ds \leq
|u(t_n)|^2 + 2\int_{t_n}^{t_{n+k}} (g, u(s)) \, ds,
\]
provided $t_{n+k} \geq t_n$.
Taking the upper limit as $k \to \infty$, we obtain
\[
\limsup_{n \to \infty} |u(t_n)|^2 + 2\nu \int_{t_n}^{t_0}\|u(s)\|^2 \, ds \leq
|u(t_n)|^2 + 2\int_{t_n}^{t_0} (g, u(s)) \, ds.
\]
Taking the lower limit as $n \to \infty$, we arrive at
\[
\limsup_{n \to \infty} |u(t_n)|^2 \leq
\liminf_{n \to \infty} |u(t_n)|^2,
\]
i.e., $\lim_{n \to \infty} |u(t_n)|$ exists. Since the limit exists for any sequence
$t_n$, it does not depend on the choice of a sequence.

For $\{t_n\}$ as in (b), the energy inequality (\ref{EI}) on
$[t_{n+k}, t_n]$ is
\[
|u(t_{n})|^2 + 2\nu \int_{t_{n+k}}^{t_n}\|u(s)\|^2 \, ds \leq
|u(t_{n+k})|^2 + 2\int_{t_{n+k}}^{t_n} (g, u(s)) \, ds,
\]
provided $t_{n+k} \leq t_n$.
Taking the lower limit as $k \to \infty$, we obtain
\[
|u(t_n)|^2 + 2\nu \int_{t_0}^{t_n}\|u(s)\|^2 \, ds \leq
\liminf_{n \to \infty} |u(t_n)|^2 + 2\int_{t_0}^{t_n} (g, u(s)) \, ds.
\]
Finally, taking the upper limit as $n \to \infty$, we arrive at
\[
\limsup_{n \to \infty} |u(t_n)|^2 \leq
\liminf_{n \to \infty} |u(t_n)|^2,
\]
i.e., $\lim_{n \to \infty} |u(t_n)|$ exists. Since the limit exists for any sequence
$t_n$, it does not depend on the choice of a sequence.
\end{proof}

\begin{lemma} \label{l:simple}
Let $u(t)$ be a Leray-Hopf solution of \eqref{NSE}
on $[T,\infty)$. Then
\[
\begin{split}
\widetilde{\lim_{t \to t_0 -}} |u(t)| =
\limsup_{t \to t_0 -} |u(t)|,\\
\widetilde{\lim_{t \to t_0 +}} |u(t)| =
\limsup_{t \to t_0 +} |u(t)|,
\end{split}
\]
and,
\begin{equation} \label{eq:orderoflim}
\widetilde{\lim_{t \to t_0 -}} |u(t)| \geq 
\widetilde{\lim_{t \to t_0 +}} |u(t)| \geq |u(t_0)|.
\end{equation}
for all $t_0 >T$.
\end{lemma}
\begin{proof}
Take any $t_0 > T$. Obviously, we have
\[
\begin{split}
\widetilde{\lim_{t \to t_0 -}} |u(t)| \leq
\limsup_{t \to t_0 -} |u(t)|,\\
\widetilde{\lim_{t \to t_0 +}} |u(t)| \leq
\limsup_{t \to t_0 +} |u(t)|.
\end{split}
\]
To show the opposite inequalities, note that for any
$t_1 \in [T,\infty) \setminus Ex$ and $t_2>t_1$,
the energy inequality (\ref{EI}) on $[t_1, t_2]$ is
\begin{equation*}
|u(t_2)|^2 + 2\nu \int_{t_1}^{t_2}\|u(s)\|^2 \, ds \leq
|u(t_1)|^2 + 2\int_{t_1}^{t_2} (g, u(s)) \, ds.
\end{equation*}
First, we fix $t_1<t_0$ and take the upper limit
as $t_2 \to t_0-$, obtaining
\[
\limsup_{t \to t_0 -} |u(t)|^2 + 2\nu \int_{t_1}^{t_0}\|u(s)\|^2 \, ds \leq
|u(t_1)|^2 + 2\int_{t_1}^{t_0} (g, u(s)) \, ds.
\]
Now we take the limit as $t_1 \to t_0 -$ avoiding the exceptional set
(see Theorem~\ref{t:lim}). We get
\[
\limsup_{t \to t_0 -} |u(t)|^2 \leq
\widetilde{\lim_{t \to t_0 -}} |u(t)|.
\]

Second, we fix $t_2>t_0$ and take the limit
as $t_1 \to t_0+$ avoiding the exceptional set. We arrive at
\[
|u(t_2)|^2 + 2\nu \int_{t_0}^{t_2}\|u(s)\|^2 \, ds \leq
\widetilde{\lim_{t \to t_0+}} |u(t)|^2 + 2\int_{t_0}^{t_2} (g, u(s)) \, ds.
\]
Taking the upper limit as $t_2 \to t_0 +$, we get
\[
\limsup_{t \to t_0 +} |u(t)|^2 \leq
\widetilde{\lim_{t \to t_0 +}} |u(t)|.
\]

Third, we fix $t_2>t_0$ and take the limit
as $t_1 \to t_0-$ avoiding the exceptional set. We obtain
\[
|u(t_2)|^2 + 2\nu \int_{t_0}^{t_2}\|u(s)\|^2 \, ds \leq
\widetilde{\lim_{t \to t_0-}} |u(t)|^2 + 2\int_{t_0}^{t_2} (g, u(s)) \, ds.
\]
Taking the limit as $t_2 \to t_0 +$ avoiding the exceptional set, we get
\[
\widetilde{\lim_{t \to t_0 +}} |u(t)|^2 \leq
\widetilde{\lim_{t \to t_0 -}} |u(t)|.
\]
Finally the weak continuity of $u(t)$ yields
\[
\widetilde{\lim_{t \to t_0 +}} |u(t)| \geq |u(t_0)|,
\]
which concludes the proof.
\end{proof}

\begin{remark} We can now rewrite the energy
inequality for a Leray-Hopf solution $u(t)$ in the following form:
\begin{equation} \label{e:gen}
|u(t)|^2 + 2\nu \int_{t_0}^{t}
\|u(s)\|^2 \, ds \leq
\widetilde{\lim_{t\to t_0+}} |u(t)|^2 + 2\int_{t_0}^{t} (g, u(s)) \, ds,
\end{equation}
for all $0\leq t_0 \leq t$.
\end{remark}

Recall that if the energy norm $|u(t)|$ of a Leray-Hopf solution 
is continuous from the right at some $t=t_0$, then $t_0$ does not belong to
the exceptional set for $u(t)$, i.e. the energy inequality holds
for $t_0$ (see Remark~\ref{r:ex}).

\begin{lemma} \label{l:contfromtheright}
Let $u(t)$ be a Leray-Hopf solution of \eqref{NSE} on $[T,\infty)$. Then
$|u(t)|$ is continuous from the right at $t=t_0 \geq T$
if and only if
\[
\widetilde{\lim_{t \to t_0 +}} |u(t)| = |u(t_0)|.
\]
\end{lemma}
\begin{proof}
If $|u(t)|$ is continuous from the right
at $t=t_0 \geq T$, then, thanks to Lemma~\ref{l:simple}, we have that
\begin{equation} \label{eq:tempcontr}
\widetilde{\lim_{t \to t_0 +}} |u(t)| = \limsup_{t \to t_0 +} |u(t)| = |u(t_0)|. 
\end{equation}
Assume now that (\ref{eq:tempcontr}) holds.
Due to the weak continuity of $u(t)$, we have
\[
\liminf_{t \to t_0 +} |u(t)|^2 \geq |u(t_0)|.
\]
Hence, 
\[
\lim_{t \to t_0 +} |u(t)|^2 = |u(t_0)|.
\]
\end{proof}

Now we will show that the strong
continuity of a Leray-Hopf solution is equivalent to the strong
continuity from the left (avoiding the exceptional set).

\begin{lemma} \label{l:simplecontin}
Let $u(t)$ be a Leray-Hopf solution of \eqref{NSE} on $[T, \infty)$. Then
$|u(t)|$ is continuous at $t=t_0>T$ if and only
if
\[
\widetilde{\lim_{t \to t_0 -}} |u(t)| = |u(t_0)|.
\]
\end{lemma}
\begin{proof}
Clearly, if $|u(t)|$ is continuous at
$t=t_0 >T$, then
\begin{equation} \label{eq:tempcont}
\widetilde{\lim_{t \to t_0 -}} |u(t)| =\limsup_{t \to t_0 -} |u(t)| =  |u(t_0)|. 
\end{equation}
Assume now that (\ref{eq:tempcont}) holds. 
Then due to the weak continuity of $u(t)$, we have
\[
\lim_{t \to t_0 -} |u(t)| = |u(t_0)|.
\]
In addition, Lemma~\ref{l:simple} (equation \eqref{eq:orderoflim}) implies that
\[
\widetilde{\lim_{t \to t_0 +}} |u(t)| = |u(t_0)|.
\]
Finally, thanks to Lemma~\ref{l:contfromtheright}, we have
\[
\lim_{t \to t_0 +} |u(t)| = |u(t_0)|.
\]
Therefore, $|u(t)|$ is continuous at $t=t_0$.
\end{proof}


We will now study a weak convergence of Leray-Hopf solutions. Our goal
is to obtain sufficient conditions for a strong convergence.

\begin{lemma} \label{l:left}
Let $\{u_n(t)\}$, $u(t)$ be Leray-Hopf solutions of \eqref{NSE} on $[T_1, \infty)$.
If $u_n \to u$ in $C([T_1, T_2]; \Hw)$, then
\[
\limsup_{n \to \infty} \widetilde{\lim_{t \to t_0 -}} |u_n(t)| \leq
\widetilde{\lim_{t \to t_0 -}} |u(t)|,
\]
for all $t_0 \in (T_1, T_2]$.
\end{lemma}
\begin{proof}
Suppose this is not true for some $t_0 \in (T_1, T_2]$. Then passing to a subsequence and dropping the subindexes, we can assume that
\[
\widetilde{\lim_{t \to t_0-}} |u_n(t)| - 
\widetilde{\lim_{t \to t_0-}} |u(t)| \geq  \delta > 0, \qquad \forall n
\]
and $|u_n(t)| \to |u(t)|$ on $[T_1, T_2]\setminus S$, where $S$ is
a set of zero measure, which includes the exceptional set for $u(t)$.

Let $S':= S\cup \left( \bigcup_n Ex_n \right)$, where 
$Ex_n$ is the exceptional set for $u_n(t)$.
The energy inequality for $u_n(t)$ implies
\[
\widetilde{\lim_{\tau \to t_0 - }} |u_n(\tau)|^2 \leq |u_n(t)|^2 +
2\int_{t}^{t_0}\left(g, u_n(s)\right) \, ds,
\]
for all $t\in [T_1,t_0] \setminus S'$.
Taking the upper limit as $n \to \infty$ and using the strong convergence
of $u_n(t)$ to $u(t)$ on $[T_1,T_2] \setminus S'$, we obtain
\[
\limsup_{n \to \infty} \widetilde{\lim_{\tau \to t_0 - }} |u_n(\tau)|^2 \leq
|u(t)|^2 + 2\int_{t}^{t_0}\left(g, u(s)\right) \, ds,
\qquad t\in [T_1, t_0] \setminus S'.
\]
Finally, letting $t \to t_0 -$, we get
\[
\limsup_{n \to \infty} \widetilde{\lim_{t \to t_0 - }} |u_n(t)|^2
\leq \widetilde{\lim_{t \to t_0-}} |u(t)|^2,
\]
which is in contradiction with the definition of $\delta$.
\end{proof}


\begin{lemma} \label{l:right}
Let $\{u_n(t)\}$, $u(t)$ be Leray-Hopf solutions of \eqref{NSE} on $[T_1, \infty)$.
If $u_n \to u$ in $C([T_1, T_2]; \Hw)$, then
\[
\liminf_{n \to \infty} \widetilde{\lim_{t \to t_0 +}} |u_n(t)| \geq
\widetilde{\lim_{t \to t_0 +}} |u(t)|,
\]
for all $t_0 \in [T_1, T_2)$.
\end{lemma}
\begin{proof}
Suppose this is not true for some $t_0 \in [T_1,T_2)$. Then passing to a
subsequence and dropping the subindexes, we can assume that
\[
\widetilde{\lim_{t \to t_0+}} |u_n(t)| - 
\widetilde{\lim_{t \to t_0+}} |u(t)| \leq \delta < 0,
\qquad \forall n
\]
and
$|u_n(t)| \to |u(t)|$ on $[T_1, T_2]\setminus S$, where $S$ is a zero
measure set, which includes the exceptional set for $u(t)$.
The energy inequality (\ref{e:gen}) for $u_n(t)$ implies
\[
|u_n(t)|^2  \leq \widetilde{\lim_{\tau \to t_0 + }} |u_n(\tau)|^2 +
2\int_{t_0}^{t}\left(g, u_n(s)\right) \, ds,
\]
for all $t\in [t_0, T_2]$.
Taking the lower limit as $n \to \infty$ and using the strong convergence
of $u_n(t)$ to $u(t)$ on $[T_1, T_2] \setminus S$, we obtain
\[
|u(t)|^2  \leq \liminf_{n \to \infty} \widetilde{\lim_{\tau \to t_0 + }}
|u_n(\tau)|^2 +
2\int_{t_0}^{t}\left(g, u(s) \right) \, ds,
\qquad t\in [t_0, T_2] \setminus S.
\]
Finally, letting $t \to t_0 +$, we get
\[
\widetilde{\lim_{t \to t_0+}} |u(t)|^2  \leq
\liminf_{n \to \infty} \widetilde{\lim_{t \to t_0 + }} |u_n(t)|^2,
\]
which is in contradiction with the definition of $\delta$.
\end{proof}

\begin{definition}
For $u(t)$, a Leray-Hopf solution of \eqref{NSE}, denote
\[
[u(t_0)] := \widetilde{\lim_{t \to t_0 -}} |u(t)| -
\widetilde{\lim_{t \to t_0 +}} |u(t)|,
\]
which we will call the energy norm jump (loss) at $t=t_0$. 
\end{definition}

Note that due to (\ref{eq:orderoflim}), the energy norm jumps of the Leray-Hopf
solutions are never negative, i.e.,
\[
[u(t)] \geq 0, \qquad \forall t,
\]
for any Leray-Hopf solution $u(t)$. Moreover, we have the following result.

\begin{theorem} \label{thm:weak-jumps}
Let $\{u_n(t)\}$, $u(t)$ be Leray-Hopf solutions of \eqref{NSE} on $[T_1, \infty)$.
If $u_n \to u$ in $C([T_1, T_2], \Hw)$, then
\[
\limsup_{n \to \infty} [u_n(t)] \leq [u(t)],
\]
for all $t \in (T_1, T_2)$.
\end{theorem}
\begin{proof}
Indeed, Lemmas \ref{l:left} and \ref{l:right}  yield
\[
\begin{split}
\limsup_{n \to \infty} [u_n(t_0)] &= \limsup_{n \to \infty}
\left(\widetilde{\lim_{t \to t_0-}}|u_n(t)| -
\widetilde{\lim_{t \to t_0+}}|u_n(t)| \right) \\
&\leq \limsup_{n \to \infty} \widetilde{\lim_{t \to t_0-}}|u_n(t)| -
\liminf_{n \to \infty} \widetilde{\lim_{t \to t_0+}}|u_n(t)| \\
&\leq \widetilde{\lim_{t \to t_0-}}|u(t)| -
\widetilde{\lim_{t \to t_0+}}|u(t)|\\
&= [u(t_0)].
\end{split}
\]
\end{proof}

Assume now that Leray-Hopf solutions converge weakly to a strongly
continuous from the right in $H$ Leray-Hopf solution. We will show that the weak 
convergence is strong if the energy jumps of solutions converge to the
energy jumps of the limit solution.

\begin{theorem} \label{t:deltat}
Let $\{u_n(t)\}$, $u(t)$ be Leray-Hopf solutions of \eqref{NSE} on $[T_1, \infty)$.
If $u_n \to u$ in $C([T_1, T_2]; \Hw)$, 
\[
\widetilde{\lim_{t \to t_0 +}}|u(t)|=|u(t_0)|, \qquad \mbox{and} \qquad
\liminf_{n \to \infty} [u_n(t_0)] \geq [u(t_0)],
\]
for some $t_0 \in (T_1, T_2)$, then
\[
\lim_{n\to \infty} |u_n(t_0)| = |u(t_0)|.
\]
\end{theorem}
\begin{proof}
Thanks to \eqref{eq:orderoflim} and Lemma~\ref{l:left}, we have
\[
\begin{split}
\limsup_{n \to \infty} |u_n(t_0)| &\leq \limsup_{n \to \infty}
\widetilde{\lim_{t\to t_0+}} |u_n(t)|\\
&\leq \liminf_{n \to \infty} \widetilde{\lim_{t\to t_0-}} |u_n(t)|
- \liminf_{n\to \infty}[u_n(t_0)]\\
&\leq\limsup_{n \to \infty} \widetilde{\lim_{t\to t_0-}} |u_n(t)| - [u(t_0)]\\
&\leq \widetilde{\lim_{t\to t_0-}} |u(t)| - [u(t_0)]\\
&= \widetilde{\lim_{t\to t_0+}} |u(t)|\\
&= |u(t_0)|,
\end{split}
\]
which concludes the proof.
\end{proof}

Note that if a Leray-Hopf solution $u(t)$ is continuous at $t=t_0$, then
$[u(t_0)]=0$. Therefore, in particular, Theorem~\ref{t:deltat}
immediately implies the following.

\begin{corollary} \label{c:gaex}
Let $\{u_n\}$ be a sequence of Leray-Hopf solutions of \eqref{NSE} on
$[T_1, \infty)$.
If $u_n \to u$ in $C([T_1, T_2]; \Hw)$, and
$|u(t)|$ is continuous at some $t=t_0 \in (T_1, T_2)$,
then $u_n(t_0) \to u(t_0)$ strongly in $H$.
\end{corollary}
So, if a solution $u(t)$ on the weak global atractor is continuous
in $H$, then it attracts its basin strongly.
Note that this is exactly Rosa's asymptotic regularity condition (see \cite{R}).

\subsection{The strong global attractor for the 3D NSE}
In this subsection we will see that if
all solutions on the weak global attractor $\Aw$ are continuous in $H$, then
$\Aw$ is the strong global attractor.
First, we will show that if a solution
belongs to $\Aw$ on some open time-interval $I$, then on any closed subinterval
of $I$ it has to coincide  with a solution that stays on $\Aw$ for all time.

\begin{lemma} \label{l:glue1}
Let $u \in \Dc([T,\infty))$ and $\tilde{u} \in \Dc([0,\infty))$, such that
$u(T)=\tilde{u}(T)$ and
\[
\widetilde{\lim_{t \to T-}} |\tilde{u}(t)| \geq 
\widetilde{\lim_{t \to T+}} |u(t)|,
\]
for some $T >0$. Let
\[
v(t)=\left\{
\begin{aligned}
&\tilde{u}(t), &0 \leq t \leq T, \\
&u(t), &t > T.
\end{aligned}
\right.
\]
Then $v \in \Dc([0,\infty))$.
\end{lemma}
\begin{proof}
Obviously, $v(t)$ is a weak solution of the 3D NSE. To show that it satisfies
the energy inequality, take any $t\geq T$ and $t_0 \in (0,T)$,
$t_0 \notin Ex$, where $Ex$ is the exceptional set for $\tilde{u}$. Then
we have
\[
|v(t)|^2 + 2\nu \int_{T}^t \|v(s)\|^2 \, ds \leq
\widetilde{\lim_{s \to T+}} |v(s)|^2 + 2\int_{T}^t (g, v(s)) \, ds,
\]
and
\[
\widetilde{\lim_{s \to T-}}|v(s)|^2 + 2
\nu \int_{t_0}^{T} \|v(s)\|^2 \, ds \leq
|v(t_0)|^2 + 2\int_{t_0}^{T} (g, v(s)) \, ds.
\]
Adding these inequalities, we obtain
\[
|v(t)|^2 + 2\nu \int_{t_0}^t \|v(s)\|^2 \, ds \leq
|v(t_0)|^2 + 2\int_{t_0}^t (g, v(s)) \, ds,
\]
which concludes the proof.
\end{proof}

\begin{corollary}
\label{l:glue}
Let $u\in \Dc([T,\infty))$ and $\tilde{u} \in \Dc((-\infty,\infty))$, such that
$u(T)=\tilde{u}(T)$ and $u(t)$ is strongly continuous from the
right at $t=T$. Let
\[
v(t)=\left\{
\begin{aligned}
&\tilde{u}(t), &t \leq T, \\
&u(t), &t > T.
\end{aligned}
\right.
\]
Then $v \in \Dc((-\infty,\infty))$.
\end{corollary}

\begin{lemma} \label{l:contona}
Let $u(t)$ be a Leray-Hopf solution $u\in \Dc ([T_1,\infty))$,
such that $u(t) \in \Aw$ for all
$t \in (T_1, T_2)$. Then for any $T_0\in(T_1,T_2)$, there
exists $v\in \Dc ((-\infty, \infty))$,
such that $u(t) = v(t)$ on $[T_0,\infty)$. In particular, $u(t) \in \Aw$ for all
$t>T_1$.
\end{lemma}
\begin{proof}
First note that $(T_1, T_0)$ contains
some interval of regularity of $u(t)$ and take $T$ in the interior of this
interval. Since $u(T) \in \Aw$, there exists a solution
$\tilde{u} \in \Dc ((-\infty, \infty))$, such that $\tilde{u}(T) =u(T)$.
We will now glue them at point $t=T$, obtaining
\[
v(t)=\left\{
\begin{aligned}
&\tilde{u}(t), &t \leq T, \\
&u(t), &t > T.
\end{aligned}
\right.
\]
Since $u(t)$ is strongly continuous at $t=T$,
Corolary~\ref{l:glue} implies that $v \in \Dc((-\infty, \infty))$.
\end{proof}

Assume now that all solutions that stay on
the weak global attractor are strongly continuous in $H$. We will prove that
in this case the weak global attractor is strong. 
This generalizes a well-known result in \cite{FT85}.
Indeed, if $\Aw$ is bounded in $V$, then all the solutions,
as long as they stay on the weak global attractor, are regular, and
are therefore continuous in $H$. This also weakens the condition
of Ball \cite{B1} -- the continuity of all Leray-Hopf solutions from $(0,\infty)$
to $H$.

\begin{lemma} \label{l:cont-comp}
If  all solutions on the weak global attractor are strongly continuous
in $H$, i.e., if $\Dc ((-\infty, \infty)) \subset C((-\infty, \infty); H)$,
then $R(t)$ is asymptotically compact.
\end{lemma} 
\begin{proof}
Take
any $\{t_n\}$, such that $t_n \to \infty$ as $n \to \infty$, and $x_n \in R(t_n)X$.
Then there exists
a sequence of solutions $v_n\in \Dc ([0,\infty))$,
such that $v_n(t_n) =x_n$.
We will show that $\{x_n\}$ has a convergent subsequence.
Without loss of generality, there exists $T>0$, such that $t_n \geq 2T$ for all $n$. Consider a sequence $u_n(t)=v_n(t+t_n -T)$, where $t\geq0$.
Due to Lemma~\ref{l:compact}, $\Dc([0,\infty))$ is compact in
$C([0, \infty);\Hw)$. Hence,
passing to a subsequence and dropping a subindex,
we can assume that $u_n$ coverges
to some $u\in \Dc ([0,\infty))$ in $C([0, \infty); \Hw)$ as $n\to \infty$.
By the definition of the weak global attractor, $u(t) \in \Aw$ for all $t \in [0, \infty)$.
Applying Lemma~\ref{l:contona} with $(T_1,T_2)=(0,T)$, we obtain that
there exists a complete trajectory
$v\in \Dc((-\infty, \infty))$, such that $u(t)=v(t)$ on $[T/2, \infty)$. Therefore,
\[
u \in C([T/2, \infty); H).
\]
Hence, Corollary~\ref{c:gaex} yields that $u_n(T) \to u(T)$ strongly in $H$,
i.e., $x_n  \to u(T)$  strongly in $H$.
\end{proof}

We can now conclude with the two main results of this section. First, as
a direct consequence of Theorem~\ref{t:asymptoticcompact} and
Lemma~\ref{l:cont-comp} we obtain

\begin{theorem} \label{thm:main}
If  all solutions on the weak global attractor are strongly continuous
in $H$, then the
strong global attractor $\As$ exists, is strongly compact, and coincides with $\Aw$.
\end{theorem} 
\begin{proof}
Due to Lemma~\ref{l:cont-comp} we have that $R(t)$ is asymptotically compact.  Then Theorem~\ref{t:asymptoticcompact} implies that
$\As$ exists, is strongly compact, and coincides with $\Aw$.
\end{proof}

Second, we prove that the condition
$\Dc((-\infty, \infty)) \subset C((-\infty, \infty);H)$
is equivalent to a condition that all the ``energy jumps'' 
uniformly converge to zero as time goes to infinity. More precisely,
let
\[
[R(t)X]:= \sup \{[u(t)] : u \in \Dc([0,\infty))\}.
\]
Then we have the following.
\begin{theorem}
$\Dc((-\infty, \infty)) \subset C((-\infty, \infty); H)$ if and only
if $[R(t)X] \to 0$ as $t \to \infty$.
\end{theorem} 
\begin{proof}
It is obvious that if  $[R(t)X] \to 0$ as $t \to \infty$, then we have
$\Dc((-\infty, \infty)) \subset C((-\infty, \infty); H)$.

Assume now that $\Dc((-\infty, \infty)) \subset C((-\infty, \infty); H)$,
but $[R(t)X]$ does not converge to $0$ as $t \to \infty$.
Then there exist a sequence of Leray-Hopf solutions
$v_n \in \Dc([0,\infty))$
and a time sequence $t_n \to \infty$ as $n\to \infty$, such that 
\begin{equation} \label{a:faraway}
\limsup_{n \to \infty}[v_n(t_n)] >0.
\end{equation}
Proceeding as in the proof of Lemma~\ref{l:cont-comp}, we can now
assume that there exists some $T>0$, such that the sequence
$u_n(t)=v_n(t+t_n-T)$ (where $t\geq 0$) converges in $C([0,\infty);\Hw)$
to the restriction $u=v|_{[0,\infty)}$ of some complete trajectory
$v \in \Dc((-\infty,\infty))$. Since $v(t)$ is continuous in $H$,
by Theorem~\ref{thm:weak-jumps} we must have
\[
\lim_{n \to \infty} [v_n(t_n)] =\lim_{n \to \infty} [u_n(T)]=[u(T)]=[v(T)]=0, 
\]
a contradiction.

\end{proof}

\subsection{Regular part of the global attractor}
We define the regular part of the weak global attractor, first introduced is
\cite{FT85}, as follows.
\[
\begin{aligned}
\Areg := \left\{ \right.&u_0:  \exists \tau>0, u \in \Dc((-\infty, \infty)) \mbox{ with }
u(0)=u_0, \mbox{ such that }
u(t) \mbox{ is regular} \\ &\mbox{on } (-\tau, \tau),
\mbox{ and for each } \tilde{u} \in \Dc((-\infty, \infty))  \mbox{ with } \tilde{u}(0)=u(0) \mbox{ we have }\\
& u(t)=\tilde{u}(t) \  \forall t \in(-\tau, \tau) \left. \right\}.
\end{aligned}
\]

The following result was proven in \cite{FT85}:
\begin{theorem} The regular part of the global attractor satisfies the
following properties:
\begin{enumerate}
\item[(a)] $\Areg$ is weakly open in $\Aw$,
\item[(b)] $\Areg$ is weakly dense in $\Aw$,
\item[(c)] If all solutions in  $\Aw$ are regular, then $\Aw$ is bounded in $V$
(hence, $\Areg=\Aw$). 
\end{enumerate}
\end{theorem}
Part $(c)$ was misstated in \cite{FT85} as
\begin{enumerate}
{\it \item[($c'$)] If $\Aw \subset V$, then $\Areg=\Aw$.}
\end{enumerate}
However, the proof provided in \cite{FT85} yields only $(c)$. As yet it is not
known whether $(c')$ is true.

Now we will further study the case when all weak solutions on the
weak global attractor of the 3D NSE are continuous in $H$. Under this
assumption, Theorem~\ref{thm:main} implies that the strong compact
global attractor $\As$ exists and $\As=\Aw$. Moreover,

\begin{theorem} If $\Dc((-\infty, \infty)) \subset C((-\infty, \infty); H)$,
then the regular part of the global attractor $\Areg$
is strongly dense in $\As$.
\end{theorem}
\begin{proof}
Due to Theorem~\ref{thm:main} $\As$ exists and is strongly compact.
Therefore, weak and strong topologies are equivalent on $\As$.
Then since $\Areg$ is weakly dense in $\Aw$, it is also strongly dense
in $\As=\Aw$.
\end{proof}

We conclude this section with the following remark. In the case when the cubic
box $\Omega=[0,L]^3$ is replaced by the cuboid $\Omega=[0,L]^2 \times [0,l]$
with $0<l \ll L$, there exists a function $\alpha(g)>0$, such that if
\begin{equation} \label{e:condonlL}
\frac{l}{L} \leq \alpha(g),
\end{equation}
then $\Aw$ is the strong global attractor and is regular (see \cite{RS,TZ}). It would be interesting to
show that for some values of $l/L$ larger than $\alpha(g)$, all the Leray-Hopf solutions
on $\Aw$ are strongly continuous, i.e., $\Aw$ is also the strong global attractor.

\section{Tridiagonal models for the Navier-Stokes equations}
In this section we introduce a two-parameter family of new simple models
for the Navier-Stokes equations with a nonlinear term enjoying the
same basic properties as the nonlinear term $B(u,u)$ in the NSE \eqref{NSE}.

The role of the space $H$ will be played by $l^2$ with the usual inner
product and norm:
\[
(u,v)= \sum_{n=1}^{\infty} u_nv_n, \qquad |u|=\sqrt{(u,u)}.
\]
The norm $|u|$ will be called the energy norm.
Let $A:D(A) \to H$ be the Laplace operator defined by
\[
(Au)_n = n^{\alpha} u_n, \qquad n\geq 1,
\]
for some $\alpha >0$. The domain $D(A)$ of this operator is
\[
\left\{u: \sum_{n=1}^\infty n^{2\alpha}u_n^2 <\infty \right\}.
\]
Clearly, $D(A)$ is dense in $H$ and 
$A$ is a positive definite operator whose eigenvalues are
\[
1, 2^{\alpha},  3^{\alpha}, \dots
\]

Let $V=A^{-1/2}H$ endowed with the following inner product and norm:
\[
((u,v))= \sum_{n=1}^{\infty} n^{\alpha}u_nv_n, \qquad \|u\|=\sqrt{((u,u))}.
\]
Here $\|u\|$ is an analog of $H^1$-norm of $u$ and we will call it the enstrophy
norm. Let also
\[
\|u\|_\gamma = \left(\sum_{n=1}^{\infty} n^\gamma u_n^2 \right)^{1/2},
\]
which is an analog of $H^{\gamma/\alpha}$-norm of $u$.

Our models for the NSE are given by the
following equations:
\begin{equation} \label{model}
\left\{
\begin{aligned}
&\ddt u_n + \nu n^{\alpha}u_n - n^\beta u_{n-1}^2 + (n+1)^\beta u_n u_{n+1}
=g_n, \qquad n=1,2,3\dots\\
&u_0=0.
\end{aligned}
\right.
\end{equation}
Here, $\nu>0$, $\alpha>0$, and $\beta>1$. Note that the value of $\ddt u_n$ is determined
only by the values of $u_{n-1}$, $u_n$, and $u_{n+1}$. Therefore, we will
refer to the equations \eqref{model} as the tridiagonal model for the Navier-Stokes
equations or shortly TNS equations.
For $u=(u_1,u_2,\dots)$ they can be written in a more condensed form as
\begin{equation}
\ddt u + \nu Au + B(u,u) = g,
\end{equation}
where
\[
(B(u,v))_n=- n^\beta u_{n-1} v_{n-1} + (n+1)^\beta u_n v_{n+1},
\]
and $u_0=0$. Note that the orthogonality property holds for $B$:
\begin{eqnarray*}
\left(B(u,v),v\right)&=&\sum_{n=1}^{\infty}\left(
-n^{\beta}u_{n-1} v_{n-1} v_n + (n+1)^\beta u_n v_{n+1} v_n \right)\\
&=&
\sum_{n=1}^{\infty}\left(
-n^{\beta}u_{n-1} v_{n-1} v_n +n^\beta u_{n-1} v_n v_{n-1} \right)\\
&=&0.
\end{eqnarray*}

In the case of TNS equations, 
a weak solution on $[T,\infty)$ (or $(-\infty, \infty)$, if
$T=-\infty$) of \eqref{model} is actually a locally bounded $H$-valued
function $u(t)$ on $[T, \infty)$, such that $u_n \in C^1([T,\infty))$
and $u_n(t)$ satisfies \eqref{model} for all $n$. From now on weak solutions will
be called just solutions.

A solution $u(t)$ is
strong (or regular) on some interval $[T_1,T_2]$, if $\|u(t)\|$ is bounded
on $[T_1,T_2]$. A solution is strong on $[T_1,\infty)$, if it is strong on every
interval $[T_1,T_2]$, $T_2\geq 0$.

A Leray-Hopf solution of \eqref{model} on the interval $[T, \infty)$
is a solution of \eqref{model} on $[T,\infty)$ satisfying the
energy inequality
\[
|u(t)|^2 + 2\nu \int_{t_0}^t \|u(\tau)\|^2 \, d\tau \leq
|u(t_0)|^2 + 2\int_{t_0}^t (g, u(\tau)) \, d\tau,
\]
for all $T \leq t_0 \leq t$, $t_0$ a.e. in $[T,\infty)$.
The set $Ex$ of those $t_0$ for which the energy
inequality does not hold will be called the exceptional set.

\subsection{A priori estimates and the existence of strong solutions}
Taking a limit of the Galerkin approximation, the existence of
Leray-Hopf solutions follows in exactly the same way as for the 3D NSE.
In this paper we will show some {\it a priori} estimates, which
can be obtained rigurously for the Leray-Hopf solutions.
For simplicity, we assume that $g$ is independent of time, $g\in H$, and
$g_n \geq 0$ for all $n$.
\\

\noindent
{\bf Energy estimates.} Formally taking a scalar product of (\ref{model}) with $u$,
we obtain
\[
\begin{split}
\frac{1}{2} \ddt |u|^2  &\leq  -\nu \|u\|^2 + |g||u|\\
&\leq  -\nu |u |^2 + \frac{\nu}{2}|u|^2 + \frac{|g|^2}{2\nu}\\
&= -\frac{\nu}{2}|u|^2 + \frac{|g|^2}{2\nu}.
\end{split}
\]
Using Gronwall's inequality, we conclude that
\begin{equation} \label{eq:defK}
|u(t)|^2 \leq e^{- \nu t} |u(0)|^2 +
\frac{|g|^2}{\nu^2}(1-e^{-\nu t}).
\end{equation}
Hence $B=\{u\in H: \ |u| \leq R\}$ is an absorbing ball for the
Leray-Hopf solutions, where $R$ is any number larger that $|g|/\nu$.
\\

\noindent
{\bf Enstrophy estimates.} Let $v=A^{1/2}u$ and 
\[
c_\mathrm{b} := \left\{
\begin{split}
&\alpha 2^\beta, &0<\alpha \leq 1,\\
&\alpha 2^{\alpha+\beta-1}, &\alpha > 1.
\end{split}
\right.
\]
Using H\"older's inequality, we obtain the following estimate
for the nonlinear term:
\[
\begin{split}
| (B(u,u),Au) | &\leq
\left| \sum_{n=1}^\infty \left[(n+1)^{\alpha}-n^{\alpha} \right]
(n+1)^{\beta }u_n^2 u_{n+1} \right|\\
&\leq c_\mathrm{b} \sum_{n=1}^\infty n^{\beta-\alpha/2-1} |v_n|^2 |v_{n+1}|\\
&\leq c_\mathrm{b} (\max_n |v_n|) \sum_{n=1}^{\infty}
n^{\beta-\alpha/2-1} v_n^2 \\
&\leq c_\mathrm{b} |v| |A^{1/2} v|^{2\beta/\alpha-2/\alpha-1} |v|^{-2\beta/\alpha+
2/\alpha+3}\\
&= c_\mathrm{b} |A^{1/2} v |^{2\beta/\alpha-2/\alpha-1}
|v|^{-2\beta/\alpha+2/\alpha+4}\\
&= c_\mathrm{b} |Au |^{2\beta/\alpha-2/\alpha-1}
\|u\|^{-2\beta/\alpha+2/\alpha+4},
\end{split}
\]
whenever $\beta\in[\alpha/2+1,3\alpha/2+1]$.
Choosing $u$ to have only two consecutive nonzero terms, it is easy to
check that this estimate is sharp. Moreover, when $\alpha=2/3$ and
$\beta=11/6$, we have
\[
| (B(u,u),Au) | \leq c_\mathrm{b} |Au |^{3/2}
\|u\|^{3/2},
\]
which is exactly what Sobolev estimates give for the 3D NSE.
Therefore, formally taking a scalar product of (\ref{model}) with $Au$,
we obtain
\[
\begin{split}
\frac{1}{2}\ddt \|u\|^2 &\leq -\nu|Au|^2 + c_\mathrm{b} |Au |^{3/2}
\|u\|^{3/2} + (g, Au)\\
&\leq -\nu|Au| + \frac{\nu}{3} |Au|^2 + \frac{3^6c_\mathrm{b}^4}{2^8\nu^3} \|u\|^6 +
\frac{3}{4\nu} |g|^2 + \frac{\nu}{3} |Au|^2\\
&\leq -\frac{\nu}{3}|Au|^2 + \frac{3^6c_\mathrm{b}^4}{2^8\nu^3} \|u\|^6 +
\frac{3}{4\nu} |g|^2,
\end{split}
\]
a Riccati-type equation for $\|u\|^2$. Hence, the model has the same
enstrophy estimates as the 3D NSE, similar properties, and the same open question concerning the regularity of the solutions in the case
$(\alpha, \beta) = (2/3, 11/6)$. In particular, we have a global existence of
Leray-Hopf weak solutions (see Theorem~\ref{thm:Leray}), local existence of
strong solutions (see Theorem~\ref{existence}), Leray's
structure theorem (see Theorem~\ref{structure}),
uniqueness of strong solutions in the class of
Leray-Hopf solutions (see Theorem~\ref{uniqueness}), and existence of
a weak global attractor (see Theorem~\ref{t:AweakNSE}).

In the case $(\alpha,\beta) = (1/2, 7/4)$, we have 
\[
| (B(u,u),Au) | \leq c_\mathrm{b} |Au |^{2}
\|u\|,
\]
which corresponds to the 4D Navier-Stokes equations.
In the case $(\alpha,\beta) = (2/5, 17/10)$, we have 
\[
| (B(u,u),Au) | \leq c_\mathrm{b} |Au |^{5/2} \|u\|^{1/2},
\]
which corresponds to the 5D Navier-Stokes equations.
In general, the choice
\[
\alpha = \frac{2}{d}, \qquad \beta = \frac{3}{2} + \frac{1}{d}
\]
would correspond to the d-dimensional Navier-Stokes equations.

The similarity of the TNS equations \eqref{model} with the NSE holds also
for values of $d\ne 3$. Indeed,
when $2\beta < 3\alpha +2$, the enstrophy estimate implies a local existence of strong solutions, i.e.,
solutions whose enstrophy norms are continuous.
More precisely, 
if $2\beta < 3\alpha +2$, then for any initial data
$u_0\in V$, there exists a
strong solution $u(t)$ with $u(0)=u_0$ on some interval $[0,T]$. 
In terms of the dimension,
a sufficient condition for the local existence of strong solutions is
$d<4.$
In the case when $\beta \leq \alpha +1$, the enstrophy estimate implies
a global regularity. In terms of the dimension,
a sufficient condition for the global existence of strong solutions is
$d\leq2$.

Now we will concentrate on the solutions with initial data $u_n(0)\geq0$
for all $n$.

\begin{theorem} \label{t:poditofsol}
Let $u(t)$ be a solution of \eqref{model} with $u_n(0)\geq 0$.
Then $u_n(t) \geq 0$ for all $t>0$, and 
$u(t)$ satisfies the energy inequality
\begin{equation} \label{ee}
|u(t)|^2 + 2\nu \int_{t_0}^t \|u(\tau)\|^2 \, d\tau \leq
|u(t_0)|^2 + 2\int_{t_0}^t (g, u(\tau)) \, d\tau
\end{equation}
for all $0 \leq t_0 \leq t$.
\end{theorem}
\begin{proof}
A general solution for $u_n(t)$ can be written as
\begin{multline*}
u_n(t)=u_n(0)\exp\left(-
\int_{0}^t\left[\nu n^\alpha + (n+1)^\beta u_{n+1}(\tau)\right] \, d\tau\right)\\
+ \int_{0}^t(g_n+n^\beta u_{n-1}^2(s) )\exp\left(-\int_{s}^{t}
\left[ \nu n^\alpha+(n+1)^\beta u_{n+1}(\tau)\right] \, d\tau \right) \, ds.
\end{multline*}
Since $u_n(0) \geq 0$ for all $n$, then  $u_n(t)\geq 0$ for all $n$, $t>0$.
Hence, multiplying \eqref{model} by $u_n$,
taking a sum from $1$ to $N$, and integrating between $t_0$ and $t$,
we obtain
\[
\begin{split}
\sum_{n=1}^N u_n(t)^2  -& \sum_{n=1}^N u_n(t_0)^2
 +2\nu \int_{t_0}^t \sum_{n=1}^N n^\alpha u_n^2 \, d\tau \\
&= - 2\int_{t_0}^t  (N+1)^\beta u_N^2 u_{N+1} \, d\tau
+2\int_{t_0}^t\sum_{n=1}^N  g_n u_n \, d\tau\\
&\leq 2\int_{t_0}^t\sum_{n=1}^N g_n u_n \, d\tau.
\end{split}
\]
Taking the limit as $N \to \infty$, we obtain \eqref{ee}.

\end{proof}

\subsection{Blow-up in finite time} \label{sub:finitetime}
Here, when $2\beta -3\alpha -3 >0$ and $g_1>0$ is large enough,
we will show that every solution $u(t)$ of \eqref{model} with $u_n(0) \geq 0$
blows up in finite time in an appropriate norm.
First, we need the following two lemmas.

\begin{lemma} \label{l:big}
Let $u(t)$ be a solution to \eqref{model} on $[0,\infty)$ with $u_n(0) \geq 0$
for all $n$. Assume that
$\|u(t)\|_{2(\beta+\gamma-1)/3} \in L^3_{\mathrm{loc}}([0,\infty);\mathbb{R})$.
Then
\begin{equation} \label{e:tmp777}
\int_{t_0}^t \sum_{n=1}^\infty n^{\beta+\gamma-1}u_n^2u_{n+1} \, d\tau < \infty,
\qquad \int_{t_0}^t \sum_{n=1}^\infty n^{\beta+\gamma-1}u_n^3
\, d\tau < \infty,
\end{equation}
and
\begin{equation} \label{e:tmp888}
\|u(t)\|_\gamma^2 - \|u(t_0)\|_\gamma^2 + 2\nu \int_{t_0}^{t} \|u\|_{\alpha+\gamma}^2 \, d\tau \geq
2\gamma \int_{t_0}^{t} \sum_{n=1}^{\infty} (n+1)^{\beta+\gamma-1}u_n^2 u_{n+1}
\, d\tau 
\end{equation}
for all $0\leq t_0 \leq t$, $0<\gamma\leq1$.
\end{lemma}
\begin{proof}
Thanks to Theorem~\ref{t:poditofsol}, $u_n(t) \geq 0$ for all $n, t>0$.
Since $\|u(t)\|^3_{2(\beta+\gamma-1)/3}$ is integrable on $[t_0,t]$  for all
$0\leq t_0\leq t$, we obtain
\begin{equation*} 
\begin{split}
\int_{t_0}^t \sum_{n=1}^\infty n^{\beta+\gamma-1}u_n^2u_{n+1} \, d\tau &\leq
2\int_{t_0}^t \sum_{n=1}^\infty n^{\beta+\gamma-1}u_n^3 \, d\tau\\
&\leq 2\int_{t_0}^t \left(\sum_{n=1}^\infty
n^{\frac{2}{3}(\beta+\gamma-1)}u_n^2\right)^{3/2} \, d\tau\\
&= 2 \int_{t_0}^{t} \|u\|_{2(\beta+\gamma-1)/3}^3 \, d\tau\\
&< \infty.
\end{split}
\end{equation*}
Hence, the relations in \eqref{e:tmp777} hold. In particular,
\begin{equation} \label{integrq}
\liminf_{n \to \infty} \int_{t_0}^t n^{\beta+\gamma} u_n^2 u_{n+1} \, d\tau =0.
\end{equation}
Now multiplying \eqref{model} by
$n^\gamma u_n$, taking a sum from $1$ to $N$, and integrating from
$t_0$ to $t$, we obtain
\begin{equation*} 
\begin{split}
\sum_{n=1}^N & n^{\gamma}u_n(t)^2  - \sum_{n=1}^N n^{\gamma} u_n(t_0)^2
+2\nu \int_{t_0}^{t} \sum_{n=1}^N n^{\alpha+\gamma} u_n^2 \, d\tau \\
&= 2\int_{t_0}^{t} \sum_{n=1}^{N-1} (n+1)^\beta((n+1)^\gamma -n^\gamma)u_n^2 u_{n+1}
\, d\tau\\
& \quad - 2\int_{t_0}^{t} (N+1)^\beta N^\gamma u_N^2 u_{N+1} \, d\tau
+ 2\int_{t_0}^{t}\sum_{n=1}^N n^\gamma g_n u_n \, d\tau\\
&\geq 2\gamma \int_{t_0}^{t} \sum_{n=1}^{N-1} (n+1)^{\beta+\gamma-1}u_n^2 u_{n+1}
\, d\tau
 - 2\int_{t_0}^{t} (N+1)^\beta N^\gamma u_N^2 u_{N+1} \, d\tau
\end{split}
\end{equation*}
Thanks to \eqref{integrq}, taking the lower limit as $N\to \infty$, we get
\eqref{e:tmp888}.

\end{proof}

\begin{lemma} \label{l:doublelarge}
For any $c>0$, there exists $g_1>0$, such that 
\[
\int_t^{t+1} |u|^2 \, d\tau > c, \qquad \forall t\geq0,
\]
for all solutions $u(t)$ with $g=(g_1,g_2,\dots)$.
\end{lemma}
\begin{proof}
Take any $c>0$. Thanks to Theorem~\ref{t:poditofsol}, $u_n(t) \geq 0$ for all $n, t>0$.
Therefore,
\[
u_1(t+1/2) \geq u_1(t) - \nu \int_t^{t+1/2} u_1 \, d\tau -
2^\beta \int_t^{t+1/2} u_1u_2 \, d\tau +g_1.
\]
Now integrating this inequality over $[t,t+1/2]$, we obtain
\[
\int_t^{t+1} u_1 \, d\tau + \nu \int_t^{t+1} u_1 \, d\tau +
2^{\beta-1} \int_t^{t+1} u_1^2 \, d\tau + 2^{\beta-1}
\int_t^{t+1} u_2^2 \, d\tau  \geq \frac{g_1}{2}. 
\]
Hence, using Cauchy-Schwarz inequality, we get
\[
(\nu+1) \left(\int_t^{t+1} (u_1^2 +u_2^2)\, d\tau \right)^{1/2} + 2^{\beta-1}
\int_t^{t+1} (u_1^2 + u_2^2) \, d\tau \geq \frac{g_1}{2}. 
\]
Obviously, for $g_1$ large enough, we have that
\[
\int_t^{t+1} |u|^2 \, d\tau \geq \int_t^{t+1} (u_1^2+u_2^2) \, d\tau > c.
\]
\end{proof}

Now we proceed to our main result in this section.

\begin{theorem} \label{thm:blowup}
For every solution $u(t)$ to equation (\ref{model}) with $u_n(0)\geq0$, $\nu>0$, $2\beta - 3\alpha -3>0$, and $g_1$ large enough,
$\|u(t)\|^3_{2(\beta+\gamma-1)/3}$ is not locally integrable for all $\gamma>0$.
\end{theorem}
\begin{proof}
Thanks to Theorem~\ref{t:poditofsol}, $u_n(t) \geq 0$ for all $n, t>0$.
Assume that there exist 
\[
\gamma\in(0,\min\{2\beta -3\alpha -3, 1\})
\]
and a solution $u(t)$ to \eqref{model}, such that
$\|u(t)\|_{2(\beta+\gamma-1)/3} \in L^3_{\mathrm{loc}}([0,\infty);\mathbb{R})$.
Then Lemma~\ref{l:big} implies that
\[
 \int_0^T \sum_{n=1}^\infty n^{\beta+\gamma-1}u_n^2u_{n+1} \, d\tau < \infty
\qquad \text{and} \qquad  \int_0^T \sum_{n=1}^\infty n^{\beta+\gamma-1}u_n^3
\, d\tau < \infty,
\]
for all $T>0$.
Moreover, since $\alpha + \gamma < 2(\beta + \gamma -1)/3$, we have that
$\|u(t)\|^2_{\alpha +\gamma}$ is locally integrable on $[0,\infty)$.

Now note that if $u_{n+1} \geq 2u_n$, then $u_nu_{n+1}^2 \leq \frac{1}{2}u_{n+1}^3$.
Otherwise, $u_nu_{n+1}^2 \leq 2 u_n^2 u_{n+1}$. Hence,
\begin{equation} \label{eq:ineq1}
u_nu_{n+1}^2 \leq {\ts \frac{1}{2}}u_{n+1}^3 + 2 u_n^2 u_{n+1}, \qquad n \in \mathbb{N}.
\end{equation}
This also implies
\begin{equation} \label{eq:ineq2}
\begin{split}
u_nu_{n+1} u_{n+2} &\leq {\ts \frac{1}{2}}u_n^2u_{n+1} + {\ts \frac{1}{2}}u_{n+1}u_{n+2}^2
\\
&\leq {\ts \frac{1}{2}}u_n^2u_{n+1} + {\ts \frac{1}{4}}u_{n+2}^3 + u_{n+1}^2u_{n+2},
\end{split}
\end{equation}
for all $n \in \mathbb{N}$. From \eqref{model} we have
\begin{multline*}
\ddt(u_nu_{n+1}) =
  -\nu (n^\alpha + (n+1)^\alpha)
u_nu_{n+1}
 +n^\beta u_{n-1}^2u_{n+1}- (n+1)^\beta u_nu_{n+1}^2\\
  +(n+1)^\beta u_n^3- (n+2)^\beta u_nu_{n+1}u_{n+2} +u_ng_{n+1}+u_{n+1}g_n.
\end{multline*}
From this, using inequalities (\ref{eq:ineq1}) and (\ref{eq:ineq2}), we obtain
\begin{multline} \label{e:ineq11}
\sum_{n=1}^\infty (n+1)^{\gamma-1}(u_nu_{n+1})(t+1)-\sum_{n=1}^\infty (n+1)^{\gamma-1}(u_nu_{n+1})(t) \\
+ 2\nu \int_t^{t+1} \sum_{n=1}^\infty (n+1)^{\alpha+\gamma-1}  u_nu_{n+1} \, d\tau
\\
+ (3+(3/2)^\beta) \int_t^{t+1}\sum_{n=1}^\infty (n+1)^{\beta+\gamma-1}  u_n^2u_{n+1} \, d\tau \\
\geq \frac{1}{4}\int_t^{t+1}\sum_{n=1}^\infty n^{\beta+\gamma -1} u_n^3 \, d\tau, 
\end{multline}
for all $t>0$.
On the other hand, Lemma~\ref{l:big} yields
\begin{equation} \label{e:ineq22}
\|u(t+1)\|^2_\gamma - \|u(t)\|^2_\gamma + 2\nu \int_t^{t+1} \|u\|_{\alpha+\gamma}^2 \, d\tau \geq
2\gamma \int_t^{t+1} \sum_{n=1}^{\infty} (n+1)^{\beta+\gamma-1}u_n^2 u_{n+1}
\, d\tau 
\end{equation}
for all $t>0$.
Denote
\[
\Theta(t) = \int_t^{t+1}\|u(\tau)\|_{\gamma}^2 \, d\tau + 
\frac{2\gamma}{3+(3/2)^\beta} \int_t^{t+1}
\sum_{n=1}^{\infty} (n+1)^{\gamma-1}(u_nu_{n+1})(\tau) \, d\tau.
\]
Note that $\Theta(t)$ is absolutely continuous on $[0,\infty)$.
We will show that $\Theta(t)$ is a Lyapunov function for the equation, i.e.,
$\Theta(t)$ is always increasing. Indeed, multiplying the inequality
(\ref{e:ineq11}) by $2\gamma/(3+(3/2)^\beta)$ and adding
(\ref{e:ineq22}), we obtain
\begin{multline*}
\ddt \Theta(t) \geq -2\nu \int_t^{t+1} \|u(\tau)\|_{\alpha + \gamma}^2 \, d\tau
- \frac{4\gamma \nu}{3+(3/2)^\beta}
\int_t^{t+1}\sum_{n=1}^{\infty} (n+1)^{\alpha+\gamma-1}  u_nu_{n+1} \, d\tau \\
+ \frac{\gamma}{6+2(3/2)^\beta}
\int_t^{t+1}\sum_{n=1}^{\infty}  n^{\beta+\gamma -1}u_n^3 \, d\tau,
\end{multline*}
a.e. on $(0,\infty)$.
Since $\gamma$ is such that
\[
\epsilon:= 2\beta -3\alpha -\gamma -3 >0,
\]
let
\[
A:=\left(\sum_{n=1}^{\infty} n^{-1-\epsilon} \right)^{-1/2}.
\]
Now H\"older's inequality yields
\[
\begin{split}
\int_t^{t+1} \sum_{n=1}^{\infty} n^{\alpha + \gamma} u_n^2 \, d\tau &\leq 
\left(\int_t^{t+1} \sum_{n=1}^{\infty} n^{-1-\epsilon} \, d\tau\right)^{1/3}
\left(\int_t^{t+1} \sum_{n=1}^{\infty} n^{\beta + \gamma -1} u_n^3 \, d\tau\right)^{2/3}
\\
&=
A^{-2/3}\left(\int_t^{t+1} \sum_{n=1}^{\infty} n^{\beta + \gamma -1} u_n^3 \, d\tau\right)^{2/3}.
\end{split}
\]
Hence,
\[
\int_t^{t+1} \sum_{n=1}^{\infty} n^{\beta + \gamma -1} u_n^3 \, d\tau \geq
A \left(\int_t^{t+1} \sum_{n=1}^{\infty} n^{\alpha + \gamma} u_n^2 \, d\tau \right)^{3/2}.
\]
Finally, we obtain
\begin{multline*}
\ddt \Theta(t)
\geq -2\nu\left( 1+ 2\gamma \frac{(3/2)^{\alpha+\gamma}}{3+(3/2)^\beta} \right)
\int_t^{t+1} \|u(\tau)\|_{\alpha + \gamma}^2 \, d\tau\\
+ \frac{\gamma A}{6+2(3/2)^\beta}
\left(\int_t^{t+1} \|u(\tau)\|_{\alpha + \gamma}^2 \, d\tau\right)^{3/2},
\end{multline*}
a.e. on $(0,\infty)$.
Due to Lemma~\ref{l:doublelarge}, if $g_1$ is large enough, then there
exists a positive constant $c$, such that
\[
\ddt \Theta(t) \geq c\Theta(t)^{3/2}, \qquad \text{a.e. on } (0,\infty).
\]
This is a Riccati-type equation. Hence, $\Theta(t)$ blows up in finite time,
which contradicts the fact that it is continuous on $[0,\infty)$.
\end{proof}

\begin{figure}
\center
\includegraphics[width=4in]{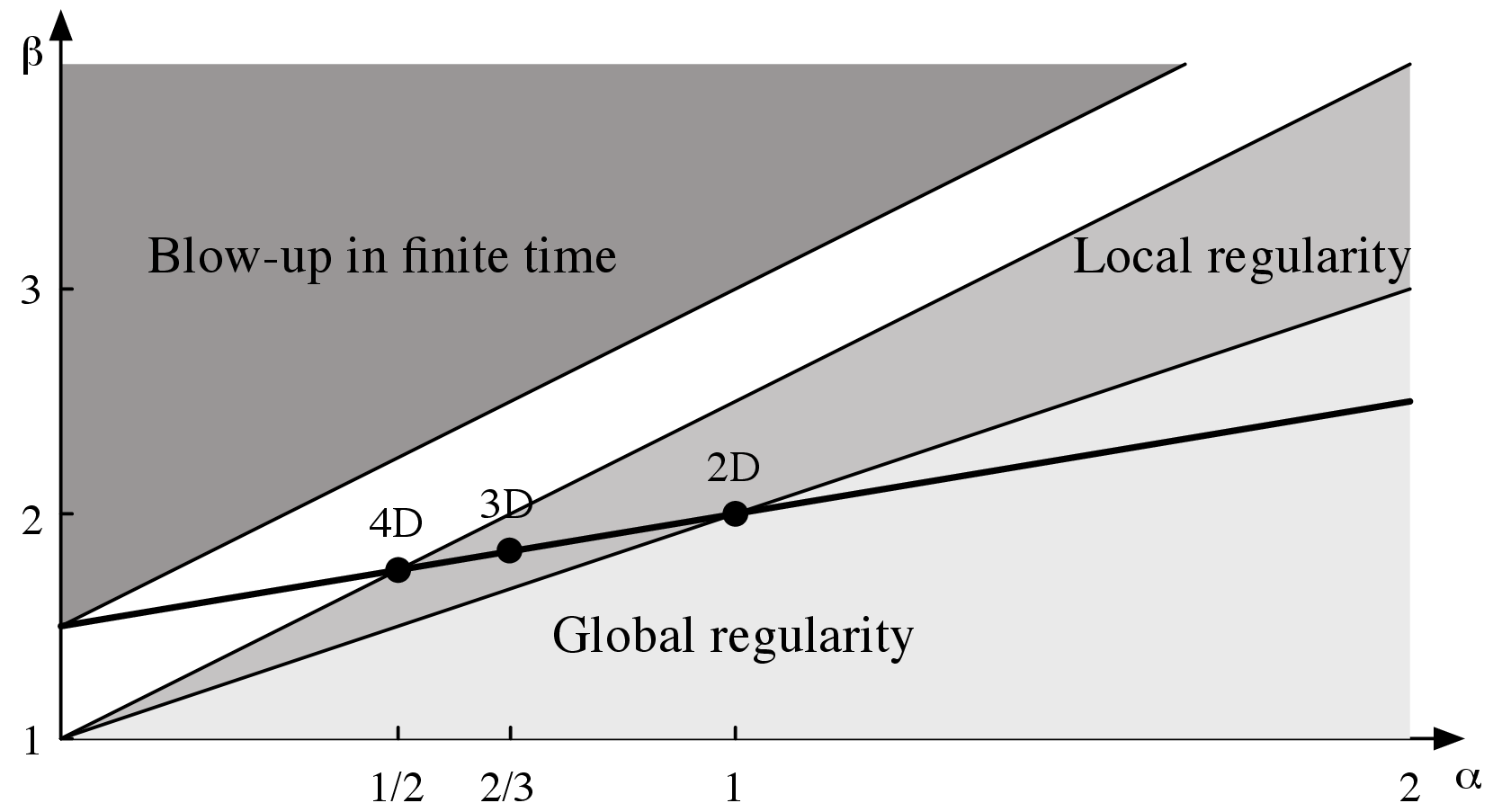}
\caption{Regions of local regularity, global regularity, and blow-up.} 
\label{fig:regions}
\end{figure}

Figure~\ref{fig:regions} shows three regions, the ones where we were able
to prove local regularity, global regularity, and  and blow-up in finite time.
The labels $2D$, $3D$, and $4D$ show the dimensions of the 
Navier-Stokes systems corresponding to the models at those points.

\subsection{Non-regular weak global attractor}
As in Subsection~\ref{s:weakattractor}, we can define an evolutionary
system $\Dc$ whose trajectories are all Leray-Hopf solutions of the
TNS equations.
The weak global attractor for this system is
\[
\begin{split}
\Aw = \{&u_0 \in H: \mbox{ there exists a Leray-Hopf solution } u(t) \mbox{ on }
(-\infty,\infty),\\
&\mbox{ such that } u(0)=u_0 \mbox{ and } |u(t)| \mbox{ is  bounded on }  (-\infty, \infty)\}.
\end{split}
\]

Recall that $g_n\geq 0$ for all $n\in \mathbb{N}$. Obviously, if $g=0$, then $\Aw=\{0\}$.
Henceforth we will assume that $g \ne 0$.

\begin{theorem} \label{t:posA}
If $g_n=0$ for all $n \geq N_g$, then
every $u=(u_1,u_2,\dots) \in \Aw$ satisfies
\[
u_n \geq 0, \qquad n=1,2,\dots
\]
\end{theorem}
\begin{proof}
A general solution for $u_n(t)$ can be written as
\begin{multline} \label{eq:gensolution}
u_n(t)=u_n(t_0)\exp\left(-
\int_{t_0}^t\nu n^\alpha + (n+1)^\beta u_{n+1}(\tau) \, d\tau\right)\\
+ \int_{t_0}^t\exp\left(-\int_{s}^{t}\nu n^\alpha+(n+1)^\beta
u_{n+1}(\tau) \, d\tau\right) (g_n+n^\beta u_{n-1}^2(s) )\, ds.
\end{multline}
Clearly, this implies the following facts.
\begin{enumerate}
\item[(a)] If $u_n(t_0) \geq 0$ for some $n$ and $t_0$, then
$u_n(t) \geq 0$ for all $t \geq t_0$.
\item[(b)] If $|u(t)|$ is bounded for all $t\in\mathbb{R}$, then
$u_n(t) \geq 0$ for all $t\in \mathbb{R}$, whenever
 $u_{n+1}(t) \geq 0$ for all $t\in \mathbb{R}$.
\end{enumerate}

Now assume that there exists  $u^0 \in \Aw$, such that $u^0_N < 0$
for some $N \geq N_g$. Then there exists a Leray-Hopf solution $u(t)$,
such that $u(0) = u^0$ and $|u(t)|$ is bounded on $(-\infty, \infty)$. For such a
solution we have $u_N(t) <0$ for all $t \leq 0$.
In addition, from the energy inequality for $u(t)$ we deduce that
\[
\begin{split}
\sum_{n=N}^{\infty} u_n(t_1)^2 - \sum_{n=N}^{\infty} u_n(t_0)^2 &\leq
2\int_{t_0}^{t_1} \left[ N^\beta u_{N-1}(\tau)^2u_{N}(\tau)
 -\nu \sum_{n=N}^\infty n^\alpha u_n(\tau)^2 \right] \, d\tau\\
 &\leq  -2 \nu \int_{t_0}^{t_1} \sum_{n=N}^\infty u_n(\tau)^2 \, d\tau,
\end{split}
\]
for all $t_0 \leq t_1 \leq 0$. Hence,
\[
\sum_{n=N}^{\infty} u_n(t_0)^2 \geq e^{-2\nu t_0}\sum_{n=N}^{\infty} u_n(0)^2,
\]
for all $t_0 \leq 0$. This implies that $|u(t)|^2$ is not bounded backwards
in time, a contradiction.

Now take any $u^0 \in \Aw$. There exists a Leray-Hopf solution $u(t)$,
such that $u(0) = u^0$ and $|u(t)|$ is bounded on $(-\infty, \infty)$.
For such a solution we proved that
\[
u_n(t) \geq 0, \qquad \forall t\in\mathbb{R}, n \geq N_g.
\]
Now note that due to the remark (b) above, if $u_{n+1}(t) \geq 0$ for all
$t\in \mathbb{R}$, then $u_n(t) \geq0$ for all $t\in \mathbb{R}$,
which concludes the proof.
\end{proof}

Now Theorem~\ref{t:posA} allows us to apply all the results in
Subsection~\ref{sub:finitetime} to the solutions on the weak global
attractor $\Aw$. Therefore, we have the following.
\begin{remark} \label{rem:bl}
Let $2\beta> 3\alpha +3$ and $g_1$ be large enough. Then
$\|u(t)\|_{2(\beta+\gamma-1)/3}$ blows up in finite time
for every solution $u(t)$ on $\Aw$, i.e., $\Aw$ is not bounded in
$H^{2(\beta + \gamma-1)/(3\alpha)}$ for any $\gamma>0$.
\end{remark}
However, this does not mean that the weak global attractor cannot be strong.
The question whether $\Aw$ is the strong global attractor remains open
in the case $\beta > \alpha +1$.

\subsection{Tridiagonal models for the Euler equations}
In this section we consider the tridiagonal models for the Euler equations
(TE), the equations \eqref{model} with $\nu=0$. First, let us show the
global existence of the weak solutions to the TE equations.
Take a sequence
$\nu_j\to 0$ as $j\to \infty$. Given $u^0\in H$,
let $u^j(t)$ be a solution of \eqref{model} with $\nu=\nu_j$ and $u^j(0)=u^0$.
It is easy
to show that the sequence $\{u^j\}$ is weakly equicontinuous. Therefore, thanks
to Ascoli-Arzela theorem, passing to a subsequence and dropping
a subindex, we obtain that there exists a function $u:[0,\infty) \to H$,
such that $u^j \to u$ in $C([0,\infty);\Hw)$
as $j\to \infty$. Clearly, $u(t)$ is a solution of the TE equations, in
the sense that it is a locally bounded $H$-valued
function on $[T, \infty)$, such that
$u_n \in C^1([0,\infty))$ and $u_n(t)$ satisfies \eqref{model}
for all $n$.

Now let us show that in the nonviscous case $\nu=0$, for every solution $u(t)$
of \eqref{model}, the norm $\|u(t)\|_{2(\beta+\gamma-1)/3}$ blows up
for any $\alpha>0$, $\beta>1$, and $\gamma>0$, reflecting the fact that there is no backward energy 
transfer for this model.

\begin{theorem} \label{t:Eblow}
Let $u(t)$ be a solution of \eqref{model} on $[0,\infty)$ with $\nu=0$,
$g_n \geq 0$, $u_n(0) \geq 0$ for all $n$, and $u(0) \ne 0$.
Then $\|u(t)\|_{2(\beta+\gamma-1)/3}$ is not bounded on $[0,\infty)$
for every $\gamma >0$ .

\end{theorem}
\begin{proof}
Clearly, it is enough to prove the theorem in the case where
$0<\gamma < \min\{1,2(\beta-1)\}$.
Assume that $\|u(t)\|_{2(\beta+\gamma-1)/3}$ is bounded on $[0,\infty)$.
Then Lemma~\ref{l:big} implies that
\begin{equation} \label{eq:dobnorme}
\|u(t)\|_\gamma^2 - \|u(t_0)\|_\gamma^2  \geq
2\gamma \int_{t_0}^{t} \sum_{n=1}^{\infty} (n+1)^{\beta+\gamma-1}u_n^2 u_{n+1}
\, d\tau \geq 0, 
\end{equation}
for all $0\leq t_0 \leq t$. Thus, $\|u(t)\|^2_{\gamma}$ is non-decreasing. 
Since $\gamma < 2(\beta+\gamma-1)/3$, $\|u(t)\|_{\gamma}$ is bounded on $[0,\infty)$. Then there
exists $E_0>0$ such that
\[
\lim_{t \to \infty} \|u(t)\|_{\gamma}^2 = E_0.
\]
Then (\ref{eq:dobnorme}) implies that 
\begin{equation} \label{eq:tempeuler}
\lim_{t\to \infty} \int_t^\infty u_n(\tau)^2u_{n+1}(\tau) \, d\tau = 0, \qquad n\in \mathbb{N}. 
\end{equation}
Hence,
\[
\begin{split}
u_n(t)^2 -u_n(0)^2&= 2n^\beta\int_0^tu_{n-1}^2u_n \, d\tau - 2(n+1)^\beta
\int_0^t u_n^2u_{n+1} \, d\tau + 2\int_0^t g_n u_n \, d\tau \\
&\to 2n^\beta\int_0^\infty u_{n-1}^2u_n \, d\tau - 2(n+1)^\beta
\int_0^\infty u_n^2u_{n+1} \, d\tau + 2\int_0^\infty g_n u_n \, d\tau,
\end{split}
\]
as $t\to \infty$. Hence $u_n(\infty) := \lim_{t\to \infty} u_n(t)$ exists for
all $n$.
Now (\ref{eq:tempeuler}) implies that $u_n(\infty) u_{n+1}(\infty)=0$ 
for all $n$. Suppose that $u_k(\infty) \ne 0$ for some $k$. Then
$u_{k+1}(\infty)=0$ and there exists $t_0>0$, such that
\[
(k+2)^{\beta}u_{k+1}(t)u_{k+2}(t) \leq  {\ts \frac{1}{3}}(k+1)^\beta u_k(\infty)^2 \qquad \mbox{and}
\qquad u_k^2(t) \geq {\ts \frac{2}{3}}u_k(\infty)^2, 
\]
for all $t\geq t_0$. Thus,
\[
\begin{split}
\ddt u_{k+1} &= (k+1)^\beta u_k(t)^2 - (k+2)^\beta u_{k+1}(t)u_{k+2}(t) +
g_{k+1} \\
&\geq {\ts \frac{1}{3}}(k+1)^\beta u_k(\infty)^2,
\end{split}
\]
for all $n \geq t_0$. Therefore,
\[
\lim_{t \to \infty} u_{k+1}(t) =\infty,
\]
a contradiction.
\end{proof}

\section*{Acknowledgment}
We thank the referee for his constructive criticism and pertinent remarks.

\begin{thebibliography}{99}


\bibitem{B1} J. M. Ball, Continuity properties and global attractors of generalized
semiflows and the Navier-Stokes equations, {\it J. Nonlinear Sci.} {\bf7} (1997),
475Ð502. Erratum: {\it J. Nonlinear Sci.} {\bf 8} (1998),  233.

\bibitem{B2} J. M. Ball, Global attractors for damped semilinear wave equations,
{\it Discr. Cont. Dyn. Sys.} {\bf10} (2004), 31Ð52.

\bibitem{CMR} T. Caraballo, P. Mar'n-Rubio, and J. C. Robinson,
A comparison between two theories for multi-valued semiflows and their
asymptotic behaviour, {\it Set-Valued Anal.} {\bf 11} (2003), 297--322.

\bibitem{CV} V. V. Chepyzhov and M. I. Vishik, {\it Attractors for Equations of
Mathematical Physics}, American Mathematical Society Colloquium Publications
{\bf 49}, American Mathematical Society, Providence, RI, 2002.

\bibitem{C} A. Cheskidov, Blow-up in finite time for the dyadic model of the
Navier-Stokes equations, {\it  Trans. Amer. Math. Soc.}, to appear.

\bibitem{CF} P. Constantin and C. Foias, {\it Navier-Stokes Equation}, University of Chicago Press, Chicago, 1989. 

\bibitem{FS}F. Flandoli and B. Schmalfu\ss, 
Weak solutions and attractors for three-dimensional Navier-Stokes equations with
nonregular force, {\it J. Dynam. Differential Equations} {\bf 11} (1999),  355--398.

\bibitem{FMRT} C. Foias, O. P. Manley, R. Rosa, and R. Temam, {\it Navier-Stokes
equatinon and Turbulence}, Encyclopedia of Mathematics and its Applications
{\bf 83}, Cambridge University Press, Cambridge, 2001. 

\bibitem{FT85}  C. Foias and R. Temam, The connection between the Navier-Stokes
equations, and turbulence theory, {\it Directions in Partial
Differential Equations} (Madison, WI, 1985), Publ. Math. Res. Center
Univ. Wisconsin, 55--73.

\bibitem{FP}
S. Friedlander and N. Pavlovi\'c, Blowup in a three-dimensional vector
model for the Euler equations, {\it Comm. Pure Appl. Math.} {\bf 57} (2004),
705--725.

\bibitem{H} J. K. Hale, {\it Asymptotic behavior of dissipative systems}, 
Amer. Math. Soc., Providence, RI, 1988.

\bibitem{HLS}J. K. Hale,  J. P. LaSalle, and M. Slemrod,
Theory of a general class of dissipative processes, {\it J. Math. Anal. Appl.}
{\bf 39} (1972), 177--191.


\bibitem{K} I. Kukavica, Role of the pressure for validity of the energy equality for solutions of
the NavierÐStokes equation,
{\it J. Dynam. Differential Equations} {\bf 18} (2006), 461--482.

\bibitem{KP}
N. H. Katz and N. Pavlovi\'c, Finite time blow-up for a dyadic model of the Euler
equations, {\it Trans. Amer. Math. Soc.} {\bf 357} (2005), 695--708.

\bibitem{L} O. A. Ladyzhenskaya, On the dynamical system generated by the Navier-Stokes
equations, {\it J. Soviet Math.}  {\bf 3} (1975), 458--479.

\bibitem{L1} O. A. Ladyzhenskaya, {\it Attractors for Semigroups and Evolution
Equations}, Cambridge Univ. Press, Cambridge, 1991.

\bibitem{MV} V. S. Melnik and J. Valero,  On attractors of multivalued semi-flows
and differential inclusions, {\it Set-Valued Anal.} {\bf 6} (1998), 83--111.

\bibitem{RS} G. Raugel and G. R. Sell, 
Navier-Stokes equations on thin 3D domains. I. Global attractors and global regularity of solutions.
{\it J. Amer. Math. Soc.} {\bf 6} (1993), 503--568

\bibitem{R} R. M. S. Rosa,  Asymptotic regularity condition for the strong
convergence towards weak limit sets and weak attractors of the 3D
Navier-Stokes equations, {\it J. Diff. Eq.}, to appear.

\bibitem{S} G. R. Sell, Global attractors for the three-dimensional Navier-Stokes
equations, {\it J. Dynam. Differential Equations} {\bf 8} (1996), 1--33. 

\bibitem{SY} G. R. Sell and Y. You, {\it Dynamics of Evolutionary Equations},
Applied Mathematical Sciences {\bf 143}, Springer-Verlag, New York, 2002.

\bibitem{T2} R. Temam, {\it Infinite Dimensional Dynamical Systems in Mechanics and
Physics}, Applied Mathematical Sciences {\bf 68}, (2nd Edition, 1997) Springer Verlag, New York, 1988.

\bibitem{TZ} R. Temam and M. Ziane, 
Navier-Stokes equations in three-dimensional thin domains with various boundary conditions.
{\it Adv. Differential Equations} {\bf 1} (1996), 499--546.
\end{thebibliography}
\end{document}